\documentclass[10pt,a4paper]{article}

\usepackage{amsmath,amssymb}
\usepackage{amsthm,latexsym}
\usepackage[labelformat=simple]{subfig}

\usepackage{bm,multirow}        
\usepackage{graphicx}
\usepackage{tcolorbox}

\graphicspath{{./fig/}}

\usepackage{hyperref}

\begin{document}

\title{Structure Preserving Model Order Reduction of Shallow Water Equations}

\author{B\"ulent Karas\"ozen\thanks{Institute of Applied Mathematics \& Department of Mathematics, Middle East Technical University, Ankara-Turkey  {bulent@metu.edu.tr}} , S\"uleyman Y{\i}ld{\i}z\thanks{Institute of Applied Mathematics, Middle East Technical University, Ankara-Turkey {yildiz.suleyman@metu.edu.tr}} , Murat Uzunca\thanks{Department of Mathematics, Sinop University, Turkey {muzunca@sinop.edu.tr}}
}

\date{}

\maketitle

\begin{abstract}

In this paper, we present two different approaches for constructing reduced-order models (ROMs) for the two-dimensional shallow water equation (SWE). The first one is based on the noncanonical Hamiltonian/Poisson form of the SWE. After integration in time by the fully implicit average vector field method, ROMs are constructed with proper orthogonal decomposition/discrete empirical interpolation method (POD/DEIM) that preserves the Hamiltonian structure. In the second approach, the SWE as a partial differential equation with quadratic nonlinearity is integrated in time by the linearly implicit Kahan's method and  ROMs are constructed with the tensorial POD that preserves the linear-quadratic structure of the SWE. We show that in both approaches, the invariants of the SWE  such as the energy, enstrophy, mass, and circulation are preserved over a long period of  time, leading to stable solutions. We conclude by demonstrating the accuracy and the computational efficiency of the reduced solutions by a numerical test problem.\\



\end{abstract}

\section{Introduction}

The shallow water equations (SWEs) consist of a set of two-dimensional partial differential equations (PDEs) describing a thin inviscid fluid layer flowing over the topography in a frame rotating about an arbitrary axis.  SWEs are widely used in modeling large-scale atmosphere/ocean dynamics and numerical weather prediction. Energy and enstrophy are the most important conserved quantities of the SWEs, whereas the energy cascades to large scales whilst enstrophy cascades to small scales \cite{Cotter18,Stewart16}. Therefore numerical schemes that preserve the energy and the enstrophy of the SWEs lead to stable solutions in the long term integration \cite{Arakawa81}. Many geophysical flows can be written in Hamiltonian form \cite{Morrison98}.
The non-canonical Hamiltonian/Poisson form of the SWE in the rotational frame with constant Coriolis force was  introduced first in \cite{Salmon04}. Later on  the Nambu formulation of the SWE \cite{Salmon07}, SWE with complete Coriolis force \cite{Salmon05a,Stewart16}, and  multi-layer SWE \cite{Stewart16} are developed.
The discrete energy conservation follows from antisymmetry
of the discrete Poisson bracket. Other conserved quantities are potential enstrophy, mass, and vorticity.

In this paper we consider two different formulations of the SWEs; as non-canonical Hamiltonian/Poisson PDE and as a PDE  with quadratic nonlinear terms in the f-plane. Both versions of the SWEs are discretized in space using finite-differences by preserving the skew-symmetry in the Poisson matrix. The system of the semi-discretized ordinary differential equation (ODE) is also in  Poisson form and has quadratic nonlinear terms. For time discretization of the SWE in Poisson form, we consider fully implicit average vector field (AVF) method  \cite{Celledoni12ped,Cohen11} that is energy
preserving. On the other hand, the ODE system with quadratic nonlinear terms obtained from the SWE in the f-plane is discretized in time by the linearly implicit Kahan's method \cite{Kahan93,Celledoni13}. Both time integrators are second-order convergent in time as the finite-difference approximation in space.  The fully implicit AVF method requires iterative solvers like Newton's method at each time step for solving the nonlinear systems arising from time discretization.   However, the linearly implicit Kahan's method requires to solve only one linear system of equations in each time step. Both methods preserve well the conserved quantities of the SWEs, like the energy, enstrophy, mass, and vorticity in long-time integration.

Numerical methods for PDEs like the SWEs are computationally expensive and require a large amount of computer memory and computing time in real-time simulations. During the last decades, reduced-order models (ROMs) have emerged as a powerful approach to reduce the cost of evaluating large systems of PDEs by constructing a low-dimensional linear subspace (reduced space), that approximately represents the solution to the system of PDEs with a significantly reduced computational cost. The solutions of the high fidelity full-order model (FOM),  generated by space-time discretization of PDEs, are projected onto reduced space using the proper orthogonal decomposition (POD). The POD has been widely used and is a computationally efficient reduced-order modeling technique in large-scale numerical simulations of nonlinear PDEs. Applying POD Galerkin projection, dominant modes of the PDEs are extracted from the snapshots of the FOM solutions. The computation of the FOM solutions is performed in the "offline" stage, whereas the reduced system from the low-dimensional subspace is solved in the "online" stage. Various ROMs are constructed for the SWEs, in conservative form  \cite{Lozovskiy16,Lozovskiy17}, in the $\beta$-plane  \cite{Stefanescu13,Stefanescu14}, in the f-plane \cite{Esfahanian09} with POD and in the $\beta$-plane \cite{Bistrian15,Bistrian17} with  dynamic mode decomposition (DMD). However, in these articles, preservation of the energy and other conservative quantities of the SWEs in the reduced space are not discussed.

The primary challenge in producing the low-dimensional models of the high-dimensional discretized PDEs is the efficient evaluation of the nonlinearities on the POD basis.
The computational cost is reduced by sampling the nonlinear terms and interpolating, known as hyper-reduction techniques \cite{Barrault04,chaturantabut10nmr,Astrid08,Zimmermann16,Nguyen08, Carlberg13}.
Our method of choice is the discrete empirical interpolation method (DEIM) \cite{chaturantabut10nmr}, which is one of the most frequently used hyper-reduction methods. The number of sampling points used in  hyper-reduction
methods often scales with the dimension of the reduced-order model, and rigorous error and stability analysis are in general not available.
For PDEs like  SWEs in the f-plane with polynomial nonlinearities,  ROMs do not require approximating the nonlinear terms through sampling. Reduced-order operators can be precomputed in the "offline" stage. Projection of FOM onto the reduced space yields low-dimensional matrix operators that preserve the polynomial structure of the FOM. This is an advantage because the offline-online computation is separated in contrast to the hyper-reduction methods.
Recently, for PDEs with polynomial nonlinearities, the computationally efficient reduced-order models are  constructed by the use of some tools from tensor theory and by matricizations of tensors \cite{Benner15,Benner18,Benner19,Gu11,Kramer19}. In this paper, we have constructed two different ROMs for  the SWEs; in Hamiltonian form  and in f-plane. Our main contribution is twofold:

\begin{itemize}

\item ROMs are constructed that preserve the reduced skew-symmetry in the Poisson matrix of the SWEs with the  energy preserving time integrator AVF method as in \cite{Gong17,Karasozen18}. We show that the reduced discrete energy (Hamiltonian) and other conserved quantities like enstrophy, mass and vorticity  are well-preserved in the long-term using POD and DEIM.

\item Semi-discrete form of the SWE in the f-plane results in an ODE with  quadratic nonlinearities, which is solved in time by linearly implicit Kahan's method. The computational complexity of the tensor products is ${\mathcal O}(n^2 N^2)$, where $n$ and $N$  are dimensions of  reduced and full order models. Applying POD in a tensorial framework (TPOD) by exploiting  matricizations of tensors, the computational cost is reduced. Consequently, the POD ROM for the quadratic system recovers an efficient offline-online decomposition and does not require an extra step of hyper-reduction, like DEIM. Here we make use of the sparse matrix technique MULTIPROD \cite{leva08mmm} to accelerate the tensor calculations.

\end{itemize}

We compare both approaches on a numerical test example. A large number of POD modes (about 50) are required for accurate ROM solutions, which is characteristic for problems with wave-type solutions like the SWEs. Numerical results demonstrate that the conserved quantities of the SWEs, energy, enstrophy, mass, and vorticity are well preserved by the FOMs and ROMs in long-time integration. Compared with the POD/DEIM ROMs with AVF discretization, the TPOD ROMs with Kahan's discretization is much faster at the same level of accuracy.

The paper organized as follows. In the next section, we briefly review the SWE in Hamiltonian form and in the f-plane. The finite-difference discretization in space and time integration by the AVF method and Kahan's method are presented in Section~\ref{fom}. POD/DEIM reduced models are constructed  for the SWE with  AVF semi-discretized Hamiltonian form, and TPOD reduced models for the semi-discrete ODE system with quadratic nonlinearities are given in Section~\ref{rom}. In Section~\ref{num}, we compare POD, POD/DEIM, and TPOD for two-dimensional rotating SWE concerning  the accuracy of the reduced solutions, preservation of the energy, enstrophy, vorticity, and computational efficiency. The paper ends with some conclusions.

\section{Rotational shallow water equation}
\label{swe}

The two-dimensional rotational SWE on a rectangular space domain $\Omega = [a,b]\times[c,d]\subset\mathbb{R}^2$ ($a,b,c,d \in \mathbb{R}$) and with the spatial element ${\mathbf x} =(x,y)^T\in\mathbb{R}^2$, is given as \cite{Salmon04}
\begin{equation} \label{swe1}
\begin{aligned}
\frac{\partial u}{\partial t } & =  qvh - \Phi_x, \\
\frac{\partial v}{\partial t } & =  -quh - \Phi_y,  \\
\frac{\partial h}{\partial t } & =   - (uh)_x   - (vh)_y,
\end{aligned}
\end{equation}
where $u({\mathbf x},t)$ and $v({\mathbf x},t)$ are the components of  the (particle) velocity in $x$ and $y$-directions, respectively, $h({\mathbf x},t)$ is the fluid depth, and
$q({\mathbf x},t) = (v_x({\mathbf x},t) - u_y({\mathbf x},t) + f({\mathbf x},t))/h({\mathbf x},t)$ is the potential vorticity with the Coriolis force $f({\mathbf x},t)$.  The subscripts $x$ and $y$ denote the partial derivatives with respect to $x$ and $y$ components, respectively. Moreover, for the gravity constant $g$, it is given that $\Phi ({\mathbf x},t) =(1/2)u({\mathbf x},t)^2 + (1/2)v({\mathbf x},t)^2 +gh({\mathbf x},t)$. The rotational SWE \eqref{swe1} is considered on a time interval $[0,T]$ for a final time $T>0$,
and with periodic boundary conditions
$$
\xi(a,y,t) = \xi(b,y,t) \; , \quad  \xi(x,c,t) = \xi(x,d,t) \; , \qquad \xi\in\{u,v,h\},
$$
and with initial conditions
$$
u({\mathbf x},0) = u_0({\mathbf x}), \quad  v({\mathbf x},0) = v_0({\mathbf x}), \quad h({\mathbf x},0) = h_0({\mathbf x}),
$$
where $u_0({\mathbf x})$, $v_0({\mathbf x})$ and $h_0({\mathbf x})$ are given functions.

\subsection{SWE in Hamiltonian form}

The non-canonical Hamiltonian/Poisson form of the SWE \eqref{swe1} is given by
\begin{equation} \label{poisson}
\dfrac{\partial z}{\partial t} = {\mathcal J}(z)\dfrac{\delta \mathcal{H}}{\delta z}
=
\begin{pmatrix}
0& q&-\partial_x \\
-q& 0&-\partial_y \\
-\partial_x&-\partial_y  & 0
\end{pmatrix}
\begin{pmatrix}
uh \\
vh \\
\frac{1}{2} (\mathbf{\upsilon}\cdot \mathbf{\upsilon})+gh
\end{pmatrix},
\end{equation}
where $z({\mathbf x},t)=(u({\mathbf x},t),v({\mathbf x},t),h({\mathbf x},t))^T$ is the solution vector and $\mathbf{\upsilon}({\mathbf x},t) = (u({\mathbf x},t),v({\mathbf x},t))^T$ is the velocity field. The (energy) functional $\mathcal{H}(z)$ denotes the
Hamiltonian given by
\begin{equation} \label{ham}
\mathcal{H}(z)= \frac{1}{2}\iint  h (\mathbf{\upsilon}\cdot \mathbf{\upsilon} +gh )  d{\mathbf x}.
\end{equation}
The skew-symmetric Poisson bracket is defined for any two functionals $\mathcal{A}$ and $\mathcal{B}$ \cite{Lynch02,Salmon04} as
\begin{equation} \label{bracket}
\{ \mathcal{A},\mathcal{B}  \} = \iint \left (q\frac{\delta((\mathcal{A},\mathcal{B})}{\delta(u,v)} - \frac{\delta \mathcal{A}}{\delta {\mathbf \upsilon}}\cdot \nabla \frac{\delta \mathcal{B}}{\delta h}
+ \frac{\delta \mathcal{B}}{\delta {\mathbf \upsilon}}\cdot \nabla \frac{\delta \mathcal{A}}{\delta h} \right ) d {\mathbf x},
\end{equation}
where $\nabla =(\partial x, \partial y)^T$, and $\delta \mathcal{A}/\delta {\mathbf \upsilon}$ is the functional derivative of $\mathcal{A}$ with respect to ${\mathbf \upsilon}$. The functional Jacobian is given by
\begin{equation*}
\frac{\delta(\mathcal{A},\mathcal{B})}{\delta(u,v)} = \frac{\delta \mathcal{A}}{\delta u}\frac{\delta \mathcal{B}}{\delta v } - \frac{\delta \mathcal{B}}{\delta u}\frac{\delta \mathcal{A}}{\delta v}.
\end{equation*}
   Although
the matrix  ${\mathcal J}$ in \eqref{poisson} is not skew-symmetric, the skew-symmetry of the Poisson
bracket appears after integrations by parts \cite{Lynch02}, and the Poisson bracket satisfies the Jacobi identity
\begin{equation*}
\{ \mathcal{A},\{\mathcal{B},\mathcal{D}\}\} + \{ \mathcal{B},\{\mathcal{D},\mathcal{A}\}\} + \{ \mathcal{A},\{\mathcal{B},\mathcal{D}\}\} = 0,
\end{equation*}
for any three functionals $\mathcal{A}$, $\mathcal{B}$ and $\mathcal{D}$.
Conservation of the Hamiltonian \eqref{ham}
follows from the antisymmetry of the Poisson bracket \eqref{bracket}
\begin{equation*}
\frac{d{\mathcal H}}{dt}= \{{\mathcal H}, {\mathcal H} \} = 0.
\end{equation*}
Other conserved quantities of the SWE \cite{Salmon04} are the Casimirs of the form
$$
\mathcal{C} = \iint  hG(q)d{\mathbf x},
$$
where $G$ is an arbitrary function of the potential vorticity
$q$.
The Casimirs are  additional constants of motion which commute with any functional $\mathcal{A}$, i.e.,  the Poisson bracket vanishes
$$
\{\mathcal{A},\mathcal{C}\} = 0, \quad \forall \mathcal{A}({\mathbf z}) \quad \text{or} \quad \mathcal{J}^{ij}\frac{\partial \mathcal{C}}{\partial {\mathbf z}^j } =0.
$$
Important special cases of the Casimirs are potential enstrophy ${\mathcal Z }$, vorticity ${\mathcal V }$ and mass ${\mathcal M }$, given by
$$
{\mathcal Z } =  \frac{1}{2}  \int \int hq^2 d{\mathbf x} =  \frac{1}{2}  \int \int \frac{1}{h}\left (\frac{\partial v}{\partial x}   - \frac{\partial u}{\partial y} + f  \right )^2 d{\mathbf x}, \quad {\mathcal V } =\iint h q d{\mathbf x}, \quad {\mathcal M } =\iint hd{\mathbf x}.
$$

\subsection{SWE in f-plane}

In almost all models of the ocean and atmosphere, for the Coriolis force, only the component of the planetary rotation vector that is locally normal to geopotential surfaces is retained. This approach is also known as "traditional approximation" \cite{Stewart16}, when the earth's surface is considered locally flat with latitudinal variation in the Coriolis force. This leads to the so-called $\beta$-plane approximation
with $f\approx f_0 \beta y$ which is used frequently in the literature  (see for example \cite{Bistrian15,Stefanescu13}). When the region is
assumed to be small enough, the latitudinal variation in the Coriolis parameter can be
ignored. Then the Coriolis force $f$ is replaced by a constant representative value, which corresponds to the f-plane formulation of the SWE.

Inserting the potential vorticity $q$ in \eqref{swe1}, the SWE can be written in the f-plane as
\begin{equation} \label{swe2}
\begin{aligned}
\frac{\partial u}{\partial t } & =  - uu_x -vu_y -gh_x + fv, \\
\frac{\partial v}{\partial t } & =  - uv_x -vv_y -gh_y -fu, \\
\frac{\partial h}{\partial t } & =  - (uh)_x -(vh)_y.
\end{aligned}
\end{equation}

\section{Full-order models}
\label{fom}

For the space discretization of the SWEs \eqref{poisson} and \eqref{swe2}, we form the uniform grid of the spatial domain $\Omega = (a,b)\times(c,d)$ with the grid nodes ${\mathbf x}_{ij} = (x_i,y_j )^T$, where $x_i=a+(i-1)\Delta x$ and $y_j=c+(j-1)\Delta y$, $i=1,\ldots, N_x+1$, $j=1,\ldots, N_y+1$, i.e., we divide the domain into $N_x$ and $N_y$ equidistant parts in $x$ and $y$-directions, respectively, with the mesh sizes $\Delta x=(b-a)/N_x$ and $\Delta y=(d-c)/N_y$. Then, we define the time-dependent semi-discrete solution vectors in the following order
\begin{equation}\label{solvec}
\begin{aligned}
{\mathbf u}(t) &= (u_{11}(t),\ldots , u_{1N_y}(t),u_{21}(t),\ldots , u_{2N_y}(t), \ldots , u_{N_xN_y}(t))^T,\\
{\mathbf v}(t) &= (v_{11}(t),\ldots , v_{1N_y}(t),v_{21}(t),\ldots , v_{2N_y}(t), \ldots , v_{N_xN_y}(t))^T,\\
{\mathbf h}(t) &= (h_{11}(t),\ldots , h_{1N_y}(t),h_{21}(t),\ldots , h_{2N_y}(t), \ldots , h_{N_xN_y}(t))^T,
\end{aligned}
\end{equation}
where each $u_{ij}(t)$, $v_{ij}(t)$ and $h_{ij}(t)$ denotes the approximation of the solutions $u({\mathbf x},t)$, $v({\mathbf x},t)$ and $h({\mathbf x},t)$, respectively, at the grid nodes ${\mathbf x}_{ij}$ at time $t$, $i=1,\ldots, N_x$, $j=1,\ldots, N_y$. Note that each semi-discrete solution vector in \eqref{solvec} has $N:=N_xN_y$ components since we omit the solutions on the most right and most top boundary because of the periodic boundary conditions. Throughout the paper, we do not explicitly represent the time dependency of the semi-discrete solutions for simplicity, and we write ${\mathbf u}$, ${\mathbf v}$ and ${\mathbf h}$, and the semi-discrete vector for the solution vector $z$ is defined by ${\mathbf z}=({\mathbf u};{\mathbf v};{\mathbf h})\in\mathbb{R}^{3N}$. Similarly, the semi-discrete vector for the potential vorticity $q$ and the Coriolis force $f$ are defined by ${\mathbf q}$ and ${\mathbf f}$, respectively.

For the approximation of the first-order partial derivative terms, we use one-dimensional central finite-differences to the first-order derivative terms in either $x$ or $y$-direction, and we extend them to two dimensions utilizing the Kronecker product.
For a positive integer $s$, let $\widetilde{D}_s$ denotes the matrix of central finite-differences  to the first-order ordinary differential operator under periodic boundary conditions
$$
 \widetilde{D}_s=
\begin{pmatrix}
 0& 1&  & &-1 \\
-1& 0&1 & &   \\
  & \ddots  & \ddots  &\ddots  &   \\
  &  &-1&0 &1 \\
 1&  &  & 1&0
\end{pmatrix} \in \mathbb{R}^{s\times s}.
$$
Then, on the two-dimensional mesh, the central finite-difference matrices corresponding to the first-order partial derivative operators $\partial_x$ and $\partial_y$ are given respectively by
$$
D_x=\frac{1}{2\Delta x}\widetilde{D}_{N_x}\otimes I_{N_y}\in\mathbb{R}^{N\times N} \; ,  \quad D_y=\frac{1}{2\Delta y}I_{N_x}\otimes \widetilde{D}_{N_y}\in\mathbb{R}^{N\times N},
$$
where $\otimes$ denotes the Kronecker product, and
 $I_{N_x}$ and $I_{N_y}$ are the identity matrices of size $N_x$ and $N_y$, respectively.

We further partition the time interval $[0,T]$ into $N_t$ uniform intervals with the step size $\Delta t=T/N_t$ as $0=t_0<t_1<\ldots <t_{N_t}=T$, and $t_k=k\Delta t$, $k=0,1,\ldots ,N_t$. Then, we denote by ${\mathbf u}^k={\mathbf u}(t_k)$, ${\mathbf v}^k={\mathbf v}(t_k)$ and ${\mathbf h}^k={\mathbf h}(t_k)$ the full discrete solution vectors at time $t_k$.

By the above formulations, the full discrete form of the conserved quantities are given as:
\begin{itemize}
\item energy
\begin{equation} \label{dener}
H^k=H({\mathbf z}^k) = \frac{1}{2}\sum_{i=1}^N (({\mathbf u}^k_i)^2  + ({\mathbf v}^k_i)^2 + g{\mathbf h}^k_i){\mathbf h}^k_i\Delta x\Delta y,
\end{equation}
\item potential enstrophy
\begin{equation} \label{denst}
Z^k=Z({\mathbf z}^k)  = \frac{1}{2} \sum_{i=1}^N \frac{((D_x{\mathbf v}^k)_i - (D_y{\mathbf u}^k)_i + {\mathbf f}^k_i)^2}{{\mathbf h}^k_i}\Delta x\Delta y,
\end{equation}
\item  vorticity and mass
\begin{equation} \label{dmv}
V^k=V({\mathbf z}^k) = \sum_{i=1}^N ((D_x{\mathbf v}^k)_i - (D_y{\mathbf u}^k)_i + {\mathbf f}^k_i)\Delta x\Delta y, \qquad M^k=M({\mathbf z}^k) = \sum_{i=1}^N {\mathbf h}^k_i\Delta x\Delta y.
\end{equation}
\end{itemize}

We can now give the FOM formulations for the SWEs \eqref{poisson} and \eqref{swe2}.

\subsection{FOM for the SWE in Hamiltonian form}

The semi-discrete formulation of the SWE \eqref{poisson} is given as the following $3N$-dimensional system of Hamiltonian ODEs
\begin{align}\label{semi-dis-poiss}
\dfrac{d {\mathbf z}}{d t}&= J({\mathbf z})\nabla H({\mathbf z})
=
\begin{pmatrix}
0& {\mathbf q}^d&-D_x \\
-{\mathbf q}^d& 0&-D_y \\
-D_x&-D_y  & 0
\end{pmatrix}
\begin{pmatrix}
{\mathbf u} \circ {\mathbf h} \\
{\mathbf v}\circ {\mathbf h} \\
\frac{1}{2} ({\mathbf u}\circ {\mathbf u} +{\mathbf v}\circ {\mathbf v})+g{\mathbf h}
\end{pmatrix},   \quad {\mathbf u},{\mathbf v},{\mathbf h} : [0,T]\mapsto \mathbb{R}^{N},
\end{align}
where $ \circ $ denotes element-wise or Hadamard product. In \eqref{semi-dis-poiss}, the matrix function ${\mathbf q}^d$ of size $N\times N$ stands for the diagonal matrix with the diagonal elements ${\mathbf q}^d_{ii}={\mathbf q}_i$ where ${\mathbf q}$ is the semi-discrete vector of the potential vorticity $q({\mathbf x},t)$, $i=1,\ldots ,N$.
Time integration of the ODE system \eqref{semi-dis-poiss} by the AVF method \cite{Cohen11} yields
\begin{align}\label{swe_avf}
 {\mathbf z}^{k+1}&={\mathbf z}^{k}+ \Delta t J\left (\dfrac{{\mathbf z}^{k+1}+{\mathbf z}^{k}}{2}\right )\int_{0}^{1} \nabla H(\xi({\mathbf z}^{k+1}-{\mathbf z}^{k})+{\mathbf z}^{k})d \xi,
\end{align}
where the integral term is computed explicitly as
\begin{align*}
&\int_{0}^{1} \nabla H(\xi({\mathbf z}^{k+1}-{\mathbf z}^{k})+{\mathbf z}^{k})d \xi\\
&=
\begin{pmatrix}
\frac{1}{3}\left({\mathbf u}^{k+1}\circ {\mathbf h}^{k+1}+ {\mathbf u}^{k}\circ {\mathbf h}^{k}\right) +
\frac{1}{6}\left({\mathbf u}^{k+1}\circ {\mathbf h}^{k}+{\mathbf u}^{k}\circ {\mathbf h}^{k+1}\right)\\
\frac{1}{3}\left({\mathbf v}^{k+1}\circ {\mathbf h}^{k+1}+{\mathbf v}^{k}\circ
{\mathbf h}^{k}\right) + \frac{1}{6}\left({\mathbf v}^{k+1}\circ {\mathbf h}^{k}+{\mathbf v}^{k}\circ {\mathbf h}^{k+1}\right)\\
\frac{1}{6}\left({\mathbf u}^{k+1}\circ {\mathbf u}^{k+1}+{\mathbf u}^{k+1}\circ {\mathbf u}^{k}+{\mathbf u}^{k}\circ {\mathbf u}^{k}\right) + \frac{1}{6}\left({\mathbf v}^{k+1}\circ {\mathbf v}^{k+1}+ {\mathbf v}^{k+1}\circ {\mathbf v}^{k}+{\mathbf v}^{k}\circ {\mathbf v}^{k}\right)
+ \frac{g}{4}\left({{\mathbf h}^{k+1}} + {{\mathbf h}^{k}}\right)
\end{pmatrix}.
\end{align*}

The AVF method \cite{Cohen11} preserves higher-order polynomial Hamiltonians, including the cubic Hamiltonian  $H$ of the SWE \eqref{swe}.
Quadratic Casimir functions like mass and circulation are preserved exactly by AVF method. But higher-order polynomial Casimirs like the enstrophy (cubic) can  not  be preserved.  Practical implementation of the AVF method requires the evaluation of the integral on the right-hand side of \eqref{swe_avf}. Since the discrete Hamiltonian $H$ and the discrete form of the Casimirs, potential enstrophy, mass and circulation are polynomial, they can be exactly integrated with a Gaussian quadrature rule of the appropriate order.
The AVF method is used with finite element discretization of the rotational SWE \cite{Cotter18,Cotter19} and for thermal SWE \cite{Eldred19} in Poisson form. Although the fully implicit integrator gives the desired properties such as conservation of the energy, it is computationally expensive. Semi-implicit implementation of the AVF  method with a simplified Jacobian with a quasi-Newton solver is used for the thermal SWE in \cite{Eldred19}. In \cite{Stewart16}, the SWE as a Poisson system is discretized in space following \cite{Salmon04} by the  Arakawa-Lamb discretization \cite{Arakawa81} in space.  It was shown that by using the time discretization with the non-structure preserving time integrators like the Adams-Bashforth and Runge-Kutta methods, drifts occur in the energy and enstrophy \cite{Stewart16}. A  general approach for constructing schemes that conserve energy and potential enstrophy formulates the SWEs  Nambu brackets \cite{Salmon05}, which is computationally expensive. Also, it is not possible to preserve multiple integrals like the enstrophy at the same time \cite{Dahlby11a} by geometric integrators like the AVF method. The SWE in a different Poisson form is discretized in space by compatible finite elements and is solved in time using the AVF method in \cite{Cotter18, Cotter19}.

\subsection{FOM for the SWE in f-plane}

Semi-discretization of the SWE \eqref{swe2} leads to the following ODE system
\begin{equation}\label{swef}
\begin{aligned}
\frac{d{\mathbf u}}{dt  }  &= - {\mathbf u}\circ (D_x{\mathbf u}) - {\mathbf v}\circ( D_y{\mathbf u})  - gD_x{\mathbf h}+{\mathbf f}\circ {\mathbf v}, \\
\frac{d{\mathbf v}}{dt }  &=  -{\mathbf u}\circ( D_x{\mathbf v})  - {\mathbf v}\circ( D_y{\mathbf v}) -gD_y{\mathbf h}-{\mathbf f}\circ {\mathbf u}, \\
\frac{d{\mathbf h}}{dt }  &=   -D_x({\mathbf u}\circ {\mathbf h})-D_y({\mathbf v}\circ {\mathbf h}).
\end{aligned}
\end{equation}
The ODE system \eqref{swef} consists of quadratic and linear parts, and can be rewritten as
\begin{align} \label{tenfom}
\frac{d{\mathbf z}}{dt} = F({\mathbf z})= R_1(\mathbf{z})+R_2(\mathbf{z})+L(\mathbf{z}),
\end{align}
where the quadratic vector fields $R_1(\mathbf{z})$ and $R_2(\mathbf{z})$ consist of the quadratic terms, while $L(\mathbf{z})$  contains the linear terms, given by
\begin{equation*}
R_1(\mathbf{z})=
\begin{pmatrix}
- {\mathbf u}\circ (D_x{\mathbf u}) \\
- {\mathbf u}\circ (D_x{\mathbf v}) \\
- D_x ({\mathbf u}\circ {\mathbf h})
\end{pmatrix} \; , \quad R_2(\mathbf{z})=
\begin{pmatrix}
- {\mathbf v}\circ (D_y{\mathbf u}) \\
- {\mathbf v}\circ (D_y{\mathbf v}) \\
- D_y ({\mathbf v}\circ {\mathbf h})
\end{pmatrix} \; , \quad L(\mathbf{z})=
\begin{pmatrix}
- gD_x{\mathbf h} + {\mathbf f}\circ {\mathbf v} \\
- gD_y{\mathbf h} - {\mathbf f}\circ {\mathbf u} \\
0
\end{pmatrix}.
\end{equation*}
For the linear-quadratic autonomous ODE systems  \eqref{tenfom},  Kahan introduced  an "unconventional"  discretization \cite{Kahan93}
\begin{equation*}
\frac{{\mathbf z}^{k+1} - {\mathbf z}^k}{\Delta t} = R_1^f({\mathbf z}^k,{\mathbf z}^{k+1})  +   R_2^f({\mathbf z}^k,{\mathbf z}^{k+1})  + \frac{1}{2}L({\mathbf z}^k + {\mathbf z}^{k+1}),
\end{equation*}
where the symmetric bilinear forms $R_1^f(\cdot ,\cdot )$  and $R_2^f(\cdot ,\cdot )$  are computed by the polarization \cite{Celledoni15}  of the quadratic vector fields $R_1(\cdot)$ and $R_2(\cdot)$, respectively, defined by
$$
R_i^f({\mathbf z}^k,{\mathbf z}^{k+1}) = \frac{1}{2}\left( R_i({\mathbf z}^k+{\mathbf z}^{k+1}) - R_i({\mathbf z}^k) -  R_i({\mathbf z}^{k+1}) \right), \quad i=1,2.
$$
Kahan's method is second-order and time-reversal and linearly implicit, i.e., it requires only one Newton iteration per time step \cite{Kahan97}
\begin{equation*}
\left ( I_{3N} -\frac{\Delta t}{2} F'({\mathbf z}^k)\right ) \frac{ {\mathbf z}^{k+1} - {\mathbf z}^k}{\Delta t}  =  F({\mathbf z}^k),
\end{equation*}
where $F'({\mathbf z}^k)$ stands for the Jacobian matrix of $F({\mathbf z})$ evaluated at ${\mathbf z}^k$, and $I_{3N}$ is the identity matrix of size $3N$. The Jacobian matrix $F'({\mathbf z})$ can be computed analytically as
\begin{equation*}
F'({\mathbf z})=-
\begin{pmatrix}
{\mathbf u}^dD_x + (D_x{\mathbf u})^d + {\mathbf v}^dD_y & (D_y{\mathbf u})^d - {\mathbf f}^d & gD_x\\
 (D_x{\mathbf v})^d + {\mathbf f}^d &  {\mathbf u}^dD_x + {\mathbf v}^dD_y + (D_y{\mathbf v})^d  & gD_y\\
 D_x {\mathbf h}^d &  D_y {\mathbf h}^d &  D_x {\mathbf u}^d + D_y {\mathbf v}^d
\end{pmatrix},
\end{equation*}
where the superscript $d$ stands for the matricization of a vector such that for any vector ${\mathbf a}\in\mathbb{R}^{N}$, the matrix ${\mathbf a}^d\in\mathbb{R}^{N\times N}$ is the diagonal matrix with the diagonal elements ${\mathbf a}^d_{ii}={\mathbf a}_i$.

A different polarization for polynomial Hamiltonians was introduced in \cite{Dahlby11} leading to linearly implicit discrete gradient methods. But this method is restricted to non-canonical Hamiltonian systems with constant Poisson matrix. Therefore it can not be applied to the rotational SWE \eqref{poisson} with a state-dependent Poisson matrix. Recently in \cite{Miyatake19} a two-step version of Kahan's method is developed for non-canonical Hamiltonians with constant Poisson matrix that preserves Hamiltonian. But this method is also not applicable to the rotational SWE due to state-dependent Poisson matrix.

\section{Reduced-order models}
\label{rom}

In this section, we give the construction of two ROMs for the rotational SWE.  The first one preserves the Hamiltonian for the SWE \eqref{poisson} in Poisson form using POD/DEIM.  The second one preserves the linear-quadratic structure of the SWE \eqref{swe2} in the f-plane representation by applying TPOD.

\subsection{Proper Orthogonal Decomposition}

Construction of ROMs is based on approximately representing the solutions of high fidelity FOMs (\eqref{semi-dis-poiss} or \eqref{swef}) from a low-dimensional linear subspace. In this study, we form the low-dimensional linear subspace from the so-called snapshot matrix of the FOM solutions applying POD, i.e., from the column space of the matrix whose columns are the discrete solution vectors at discrete time instances.
In reduced-order modeling of fluid dynamics problems, to avoid the fluctuations,  the mean centered
snapshots are collected at time instances $t_k$, $k=1,\cdots, N_t$, in the snapshot matrices $S_u$, $S_v$ and $S_h$
\begin{align*}
S_u &= \left({\mathbf u}^1 - \overline{\mathbf u}, {\mathbf u}^2 - \overline{\mathbf u}, \cdots, {\mathbf u}^{N_t} - \overline{\mathbf u} \right)\in\mathbb{R}^{N\times N_t} ,\\
S_v &= \left({\mathbf v}^1 - \overline{\mathbf v}, {\mathbf v}^2 - \overline{\mathbf v}, \cdots, {\mathbf v}^{N_t} - \overline{\mathbf v}\right)\in\mathbb{R}^{N\times N_t} ,\\
S_h &= \left({\mathbf h}^1 - \overline{\mathbf h}, {\mathbf h}^2 - \overline{\mathbf h}, \cdots, {\mathbf h}^{N_t} - \overline{\mathbf h}\right)\in\mathbb{R}^{N\times N_t},
\end{align*}
where $\overline{\mathbf{u}}$, $\overline{\mathbf{v}}$, $\overline{\mathbf{h}}\in\mathbb{R}^{N}$  denote the mean of the snapshots defined by
$$
\overline{\mathbf u} = \frac{1}{N_t}\sum_{k=1}^{N_t} { \mathbf u}^k\; , \quad \overline{\mathbf v} = \frac{1}{N_t}\sum_{k=1}^{N_t} { \mathbf v}^k\; , \quad \overline{\mathbf h} = \frac{1}{N_t}\sum_{k=1}^{N_t} { \mathbf h}^k.
$$
For a positive integer $n\ll N_t$, the POD aims to form the $n$-POD basis matrices $V_{u,n}, V_{v,n}, V_{h,n}\in\mathbb{R}^{N\times n}$ which approximately span the column space of the snapshot matrices $S_u$, $S_v$ and $S_h$, respectively.
The $n$-POD basis matrices are computed through the application of the singular value decomposition (SVD) to the snapshot matrices
\begin{equation*}
S_u= V_u \Sigma_u W_u^T\; , \quad S_v= V_v \Sigma_v W_v^T \; , \quad S_h= V_h \Sigma_h W_h^T \; ,
\end{equation*}
where the columns of the orthonormal matrices $V_i \in \mathbb{R}^{ N\times N_t}$ and
$W_i\in \mathbb{R}^{ N_t\times N_t}$ are the left and right singular vectors of the snapshot matrix $S_i$, respectively, and the diagonal matrix $\Sigma_i\in \mathbb{R}^{ N_t\times N_t}$ with the diagonal elements $(\Sigma_i)_{jj}= \sigma_{i,j}$, $j=1,\ldots ,N_t$, contains the singular values of $S_i$, $i\in\{u,v,h\}$. Then, the $n$-POD basis matrix $V_{i,n}$ consists of the first $n$ left singular vectors (POD modes) from $V_i$, $i\in\{u,v,h\}$.
The $n$-POD basis matrix minimizes the following least squares error
\begin{equation*}
\min_{V_{i,n}\in \mathbb{R}^{ N\times n}}||S_i-V_{i,n}V_{i,n}^TS_i ||_F^2 = \sum_{j=n+1}^{N_t} \sigma_{i,j}^2\; , \quad i\in\{u,v,h\},
\end{equation*}
where $\|\cdot\|_F$ denotes the Frobenius norm. Thus the error in the snapshot representation is given by the sum of the squares of the singular values corresponding to those left singular vectors which are not included in the $n$-POD basis matrix. The singular values provide quantitative guidance for choosing the size of the POD basis, that accurately represents the given snapshot data. Usually, the following relative "cumulative energy" criterion is used
\begin{equation} \label{energy_criteria}
\frac{\sum_{j=1}^n \sigma_{i,j}^2}{\sum_{j=1}^{N_t} \sigma_{i,j}^2  } > 1 - \kappa,
\end{equation}
where $\kappa $ is a user-specified tolerance, typically taken to be $\kappa\leq 10^{-3}$ to catch at least $99.9 \%$ of data information.

After computing the $n$-POD basis matrices $V_{u,n}$, $V_{v,n}$ and $V_{h,n}$, we can construct reduced approximations $\widehat{\mathbf u}$, $\widehat{\mathbf v}$ and $\widehat{\mathbf h}$ to the FOM solutions ${\mathbf u}$, ${\mathbf v}$ and ${\mathbf h}$, respectively, which are given by
\begin{equation*}
{\mathbf u} \approx \widehat{\mathbf u}=\overline{\mathbf u}+ V_{u,n}  {\mathbf u}_r \; , \quad {\mathbf v} \approx \widehat{\mathbf v}=\overline{\mathbf v}+ V_{v,n}  {\mathbf v}_r \; , \quad {\mathbf h} \approx \widehat{\mathbf h}=\overline{\mathbf h}+ V_{h,n}  {\mathbf h}_r \; ,
\end{equation*}
where  ${\mathbf u}_r$, ${\mathbf v}_r$, ${\mathbf h}_r\in\mathbb{R}^{n}$ are the solution vectors of the reduced system of dimension $n\ll N$.
For convenience, we also define the following vectors and matrix
\begin{equation*}
\overline{\mathbf z} =
\begin{pmatrix}
\overline{\mathbf u}\\
\overline{\mathbf v}\\
\overline{\mathbf h}
\end{pmatrix}\in\mathbb{R}^{3N} \; , \quad
{\mathbf z}_r =
\begin{pmatrix}
{\mathbf u}_r\\
{\mathbf v}_r\\
{\mathbf h}_r
\end{pmatrix}\in\mathbb{R}^{3n} \; , \quad V_{z,n}=
\begin{pmatrix}
V_{u,n} & &\\
& V_{v,n}&\\
& & V_{h,n}
\end{pmatrix}\in\mathbb{R}^{3N\times 3n},
\end{equation*}
and then we get the approximation ${\mathbf z} \approx \widehat{\mathbf z}=\overline{\mathbf z}+ V_{z,n}  {\mathbf z}_r$. We also note that since the columns of the $n$-POD matrices are orthonormal, the reverse of the approximation identity holds as ${\mathbf z}_r = V_{z,n}^T(\widehat{\mathbf z} -\overline{\mathbf z})$.

\subsection{ROM for the SWE in Hamiltonian form}

The semi-discretized SWE \eqref{semi-dis-poiss} in Hamiltonian form can be written as a nonlinear ODE of the form
\begin{equation} \label{fomham}
\dfrac{d {\mathbf z}}{d t} = F({\mathbf z}) \; , \qquad F({\mathbf z})= \begin{pmatrix}
F_1({\mathbf z}) \\
F_2({\mathbf z}) \\
F_3({\mathbf z})
\end{pmatrix}=
\begin{pmatrix}
{\mathbf q}\circ{\mathbf v}\circ {\mathbf h} - D_x \left( \frac{1}{2}({\mathbf u}\circ {\mathbf u} +{\mathbf v}\circ {\mathbf v})+g{\mathbf h} \right)     \\
-{\mathbf q}\circ{\mathbf u}\circ {\mathbf h} - D_y\left( \frac{1}{2}({\mathbf u}\circ {\mathbf u} +{\mathbf v}\circ {\mathbf v})+g{\mathbf h}\right)  \\
-D_x({\mathbf u}\circ {\mathbf h}) - D_y({\mathbf v}\circ {\mathbf h})
\end{pmatrix}.
\end{equation}
By substituting the identity ${\mathbf z} \approx \widehat{\mathbf z}=\overline{\mathbf z}+ V_{z,n}  {\mathbf z}_r$ into  \eqref{fomham}, applying the  Galerkin
projection onto $V_{z,n}$, and using the orthonormality fact $V_{z,n}^TV_{z,n}=I_{3n}$, the following $3n$-dimensional POD reduced model is obtained
\begin{equation}\label{rom-semi}
\frac{d {\mathbf z}_r}{dt } = V_{z,n}^T F(\widehat{\mathbf z})=\begin{pmatrix}
V_{u,n}^TF_1(\widehat{\mathbf z}) \\
V_{v,n}^TF_2(\widehat{\mathbf z}) \\
V_{h,n}^TF_3(\widehat{\mathbf z})
\end{pmatrix}, \qquad {\mathbf z}_r(0) = V_{z,n}^T({\mathbf z}(0) -\overline{\mathbf z}).
\end{equation}

Although the dimension of the ROM \eqref{rom-semi} is $3n\ll 3N$, the cost for evaluating the nonlinear vector $F(\cdot)$ scales not only with the reduced dimension $n$ but also with the dimension $N$ of the FOM. The computational cost is reduced by sampling the nonlinearity $F(\cdot)$ and then interpolating, known as hyper-reduction technique. Several hyper-reduction methods are developed to reduce the computational cost of evaluating the reduced nonlinear terms at selected points.
We use DEIM \cite{chaturantabut10nmr} which is the most frequently used hyper-reduction methods.
By the DEIM, full discrete nonlinear vectors are collected in the snapshot matrices defined by
$$
G_i = ( F_i^1,F_i^2,\cdots , F_i^{N_t})\in\mathbb{R}^{N\times N_t}, \quad i=1,2,3,
$$
where $F_i^k\in\mathbb{R}^N$ denotes the $i$-th component of the nonlinear vector $F({\mathbf z})$ in \eqref{fomham} at time $t_k$, $k=1,\ldots , N_t$. Then, we can approximate each $F_i({\mathbf z})$ in the column space of the snapshot matrix $G_i$. We first apply POD to the snapshot matrix $G_i$ and find the basis matrix $V_{F_i,m}\in\mathbb{R}^{N\times m}$ whose columns are the basis vectors (DEIM modes) spanning the column space of the snapshot matrix $G_i$.  Then, we apply the DEIM algorithm \cite{chaturantabut10nmr} to find a projection matrix $P_i\in\mathbb{R}^{N\times m}$ so that we have the approximation
$$
F_i({\mathbf z}) \approx V_{F_i,m}(P_i^TV_{F_i,m})^{-1} P_i^TF_i({\mathbf z}),
$$
and then we get the DEIM approximation to the reduced nonlinearities in \eqref{rom-semi} as
$$
V_{u,n}^TF_1(\widehat{\mathbf z}) \approx \mathcal{V}_{u,1} (P_1^TF_1(\widehat{\mathbf z})) \; , \quad V_{v,n}^TF_2(\widehat{\mathbf z}) \approx \mathcal{V}_{v,2} (P_2^TF_2(\widehat{\mathbf z})) \; , \quad V_{h,n}^TF_3(\widehat{\mathbf z}) \approx \mathcal{V}_{h,3} (P_3^TF_3(\widehat{\mathbf z})),
$$
where $\mathcal{V}_{u,1} = V_{u,n}^TV_{F_1,m}(P_1^TV_{F_1,m})^{-1}$, $\mathcal{V}_{v,2} = V_{v,n}^TV_{F_2,m}(P_2^TV_{F_2,m})^{-1}$ and $\mathcal{V}_{h,3} = V_{h,n}^TV_{F_3,m}(P_3^TV_{F_3,m})^{-1}$ are all the matrices of size $n\times m$, and they are precomputed in the offline stage. Note that the terms $P_i^TF_i(\widehat{\mathbf z})\in\mathbb{R}^{m}$ means that we need to compute only $m\ll N$ entry of the nonlinear vector $F_i(\widehat{\mathbf z})$, $i=1,2,3$. In addition, the computational complexity for the Jacobian matrix reduces from ${\mathcal O} (N^2)$ to ${\mathcal O} (nN)$. By using DEIM approximation, the semi-discrete ROM \eqref{rom-semi} with POD/DEIM becomes
\begin{equation}\label{rom-semi-deim}
\frac{d {\mathbf z}_r}{dt } = \begin{pmatrix}
\mathcal{V}_{u,1}P_1^TF_1(\widehat{\mathbf z}) \\
\mathcal{V}_{v,2}P_2^TF_2(\widehat{\mathbf z}) \\
\mathcal{V}_{h,3}P_3^TF_3(\widehat{\mathbf z})
\end{pmatrix}, \qquad {\mathbf z}_r(0) = V_{z,n}^T({\mathbf z}(0) -\overline{\mathbf z}).
\end{equation}
The ROM \eqref{rom-semi-deim} is also integrated in time by the AVF method similar to the scheme used for the FOM \eqref{semi-dis-poiss}.

In the case of the selection of the number $m$ of DEIM modes, the "cumulative energy" criterion \eqref{energy_criteria} is used. But, because of the nature of the nonlinearity, snapshot matrices are more sensitive and a larger number of modes are needed for accurate approximation. We take $\kappa\leq 10^{-4}$ to catch at least $99.99 \%$ of data.

\subsection{ROM for the SWE in f-plane}

For PDEs and ODEs with polynomial nonlinearity, ROMs do not require approximating the nonlinear function $F(\cdot)$ through sampling, the reduced-order operators can be precomputed in the offline stage. This is beneficial because the offline-online computation is separated in contrast to the hyper-reduction methods. In the past, for the Navier-Stokes \cite{Graham99,Holmes12} equation, the quadratic polynomial forms of the FOMs are exploited by constructing reduced models.  For polynomial nonlinearity, this avoids the approximation of the nonlinear terms by hyper-reduction and allows separation of offline and online computation of FOM and ROM, and the reduced model preserves the linear-quadratic structure of the original system.

In order to avoid the approximation of the nonlinear terms in the SWE \eqref{tenfom} by hyper-reduction, such as  DEIM, we represent the system in terms of tensors and Kronecker product.
On the other hand, mathematical operations with tensors can be easily performed using their corresponding matrix
representations. Unfolding a tensor into a matrix is called matricization of a tensor.  A common matricization of a tensor ${\mathcal Q}$ is the $\mu$-mode matricization $Q^{(\mu)}$ \cite{Benner15}. For a third-order tensor in quadratic systems like \eqref{tenfom}, there
are three different ways of unfolding, depending on the $\mu$-mode that are used for the unfolding. For example $Q^{(1)}$ is called $1$-mode matricization of ${\mathcal Q}$. $2$-mode and $3$-mode matricizations can be constructed similarly \cite{Benner15}. We refer to \cite{Kolda09} for more details on these basic concepts of
tensors. Another advantage of matricizations is that tensor-matrix
multiplications can be performed by matrix-matrix products.
The special structure of the matrix $Q$ which represents the Hessian of the right-hand side of \eqref{tenfom} has been exploited numerically in an  efficient way to construct the reduced-order Hessian \cite{Benner15}. The Hessian matrix $Q$ is an unfolding of a $3$-tensor ${\mathcal Q}\in \mathbb{R}^{3N\times 3N\times 3N}$. It is sparse due to the local structure of common discretization like finite-differences, finite elements for polynomial nonlinearities. For the SWE with quadratic nonlinear terms, the cross terms $\mathbf{z}_i\cdot \mathbf{z}_j$ vanish for $|i-j|>3$ in the semi-discretized ODE \eqref{tenfom}.  Therefore the number of nonzero columns of $Q$ is only $2N$.

The semi-discrete system \eqref{tenfom} with quadratic and linear parts can be rewritten using the Kronecker product as the following
\begin{align}\label{tenfom2}
\frac{d{\mathbf z}}{dt} = F({\mathbf z})= \widetilde{R}_1(\mathbf{z})+\widetilde{R}_2(\mathbf{z})+L(\mathbf{z}),
\end{align}
with the following redefined quadratic parts
\begin{equation}\label{kron}
\widetilde{R}_1(\mathbf{z})=-A^xQ\left(
\begin{pmatrix}
{\mathbf u}\\
{\mathbf u}\\
{\mathbf u}\\
\end{pmatrix}
\otimes(B^x\mathbf{z})\right), \qquad \widetilde{R}_2(\mathbf{z})=-A^yQ\left(
\begin{pmatrix}
{\mathbf v}\\
{\mathbf v}\\
{\mathbf v}\\
\end{pmatrix}
\otimes(B^y\mathbf{z})\right),
\end{equation}
where $\otimes$ denotes the Kronecker product, and the matrices $A^x, A^y, B^x, B^y\in\mathbb{R}^{3N\times 3N}$ are given by
\begin{align*}
& A^x=
\begin{pmatrix}
I_N & &\\
& I_N&\\
& & D_x
\end{pmatrix},\quad
A^y=
\begin{pmatrix}
I_N & &\\
& I_N&\\
& & D_y
\end{pmatrix},\\
& B^x=
\begin{pmatrix}
D_x & &\\
&D_x&\\
& & I_N
\end{pmatrix},\quad
B^y=
\begin{pmatrix}
D_y & &\\
&D_y&\\
& & I_N
\end{pmatrix},
\end{align*}
where $I_N$ denotes the identity matrix of size $N\times N$. In \eqref{kron}, the matrix  $Q \in \mathbb{R}^{3N \times (3N)^2}$ represents the matricized $3$-tensor such that $Q(\mathbf{z}\otimes \mathbf{z})=\mathbf{z}\circ \mathbf{z}$ is satisfied.
By substituting the identity ${\mathbf z} \approx \widehat{\mathbf z}=\overline{\mathbf z}+ V_{z,n}  {\mathbf z}_r$ into \eqref{tenfom2}, and applying Galerkin projection onto $V_{z,n}$, we obtain the following reduced linear-quadratic equation
\begin{equation}\label{quad1}
\frac{d {\mathbf z}_r}{dt} = F_r(\widehat{\mathbf z}) = F_r^u(\widehat{\mathbf z})  + F_r^v(\widehat{\mathbf z}) + L_r(\widehat{\mathbf z}),
\end{equation}
where we set
$$
F_r^u(\widehat{\mathbf z}) = -V_{z,n}^TA^xQ\left(
\begin{pmatrix}
\widehat{\mathbf u}\\
\widehat{\mathbf u}\\
\widehat{\mathbf u}\\
\end{pmatrix}
\otimes(B^x\widehat{\mathbf z})\right) \; , \quad F_r^v(\widehat{\mathbf z}) = -V_{z,n}^TA^yQ\left(
\begin{pmatrix}
\widehat{\mathbf v}\\
\widehat{\mathbf v}\\
\widehat{\mathbf v}\\
\end{pmatrix}
\otimes(B^y\widehat{\mathbf z})\right)\; , \quad L_r(\widehat{\mathbf z})= V_{z,n}^TL(\widehat{\mathbf z}).
$$
In the reduced linear-quadratic equation \eqref{quad1}, the tricky part is the computation of the reduced quadratic parts $F_r^u(\widehat{\mathbf z})$ and $F_r^v(\widehat{\mathbf z})$. We next discuss the computation of the reduced quadratic part $F_r^u(\widehat{\mathbf z})$, and then the other reduced quadratic part $F_r^v(\widehat{\mathbf z})$ can be computed in a similar way. Substituting the identity $\widehat{\mathbf z}=\overline{\mathbf z}+V_{z,n}{\mathbf z}_r$ into the term $F_r^u(\widehat{\mathbf z})$, and using the properties of the Kronecker product operation, we obtain
\begin{equation} \label{tencomp}
\begin{aligned}
F_r^u(\widehat{\mathbf z})= F_r^u(\overline{\mathbf z}+V_{z,n}{\mathbf z}_r) &= -V_{z,n}^TA^xQ\left(
\begin{pmatrix}
\overline{\mathbf{u}}+V_{u,n} \mathbf{u}_r\\
\overline{\mathbf{u}}+V_{u,n} \mathbf{u}_r\\
\overline{\mathbf{u}}+V_{u,n} \mathbf{u}_r\\
\end{pmatrix}
\otimes (B^x(\overline{\mathbf{z}}+V_{z,n} \mathbf{z}_r))\right), \\
&=- V_{z,n}^TA^xQ
\left(
\begin{pmatrix}
\overline{\mathbf{u}}\\
\overline{\mathbf{u}}\\
\overline{\mathbf{u}}\\
\end{pmatrix}
\otimes(B^x\overline{\mathbf{z}})+
\begin{pmatrix}
\overline{\mathbf{u}}\\
\overline{\mathbf{u}}\\
\overline{\mathbf{u}}\\
\end{pmatrix}
\otimes(B^xV_{z,n} \mathbf{z}_r) \right) \\
& \quad - V_{z,n}^TA^xQ
	\left(\begin{pmatrix}
	V_{u,n} \mathbf{u}_r\\
	V_{u,n} \mathbf{u}_r\\
	V_{u,n} \mathbf{u}_r\\
	\end{pmatrix}
	\otimes(B^x\overline{\mathbf{z}})+
	\begin{pmatrix}
	V_{u,n} \mathbf{u}_r\\
	V_{u,n} \mathbf{u}_r\\
	V_{u,n} \mathbf{u}_r\\
	\end{pmatrix}
	\otimes(B^xV_{z,n} \mathbf{z}_r) \right).
	\end{aligned}
\end{equation}
In \eqref{tencomp}, only the last term in the last line is quadratic and the other terms are at most linear. The quadratic term is computed by the Kronecker product as follows
\begin{equation} \label{tencomp2}
	\begin{aligned}
	-V_{z,n}^TA^xQ
	\left(
	\begin{pmatrix}
	V_{u,n} \mathbf{u}_r\\
	V_{u,n} \mathbf{u}_r\\
	V_{u,n} \mathbf{u}_r\\
	\end{pmatrix}
	\otimes(B^xV_{z,n} \mathbf{z}_r) \right) &=-V_{z,n}^TA^xQ\left(
	V_{u,n}^*\otimes (B^xV_{z,n}) \right)
	\left(
	\begin{pmatrix}
	\mathbf{u}_r\\
	\mathbf{u}_r\\
	\mathbf{u}_r\\
	\end{pmatrix}
	\otimes( \mathbf{z}_r) \right)\\
	&=Q_{u,r}\left(
	\begin{pmatrix}
	\mathbf{u}_r\\
	\mathbf{u}_r\\
	\mathbf{u}_r\\
	\end{pmatrix}
	\otimes( \mathbf{z}_r) \right),
	\end{aligned}
\end{equation}
where the matrix $Q_{u,r}=-V_{z,n}^TA^xQ\left(
	V_{u,n}^*\otimes (B^xV_{z,n}) \right)\in\mathbb{R}^{3n\times (3n)^2}$ can be computed in the offline stage, and
	\begin{align*}
	V_{u,n}^*=
	\begin{pmatrix}
	V_{u,n}& &\\
	& V_{u,n}&\\
	& & V_{u,n}
	\end{pmatrix}\in\mathbb{R}^{3N\times 3n}.
	\end{align*}

The main computational burden is computation of the Kronecker product $V_{u,n}^*\otimes (B^xV_{z,n})\in \mathbb{R}^{(3N)^2\times (3n)^2}$ in $Q_{u,r}$, which has complexity of order ${\mathcal O}(n^2N^2)$ for quadratic nonlinearity, and computation of  $Q_{u,r}$  in \eqref{tencomp2} is costly  due to the dense structure of the POD basis matrices.  In \cite{Benner15}, an algorithm  is developed to construct the reduced matricized tensor $Q_{u,r}$ for quadratic nonlinearity which avoids the computation of the Kronecker product $V_{u,n}^*\otimes (B^xV_{z,n})$, having a complexity of order ${\mathcal O}(nN^2)$.
Therein, using the $\mu$-mode (matrix) product the reduced matricized tensor can be efficiently computed as following
	\begin{itemize}
		\item Compute $ \mathcal{Y}^{3n\times 3N \times 3N} $ by $ Y^{(1)}= -V_{z,n}^TA^xQ $,
		\item Compute $ \mathcal{Z}^{3n\times 3n\times 3N} $ by $ Z^{(2)}= V_{z,n}^T(B^x)^TY^{(2)}$,
		\item Compute $ \mathcal{Q}_{u,r}^{3n\times 3n\times 3n} $ by $ Q_{u,r}^{(3)}=(V_{u,n}^*)^T Z^{(3)}$.
	\end{itemize}

Although the $ \mu $-mode (matrix) product decreases the complexity of evaluating the reduced matrix $ {Q}_{u,r} $, still the matrix ${Q} $ has to be built for each different polynomial nonlinearity.
Recently two new algorithms \cite{Benner18,Benner19} are developed for more efficient computation of the reduced matrix $ {Q}_{u,r} $  using  the particular structure of Kronecker product.
 A more compact form of the evaluation of the reduced matrix $ {Q}_{u,r} $ is given in MatLab notation as follows  \cite{Benner18,Benner19}
\begin{equation}\label{goyal}
	\begin{aligned}
	{Q}_{u,r} =-V_{z,n}^T A^xQ
	\left(V_{u,n}^* \otimes (B^xV_{z,n}) \right)
	=-V_{z,n}^T A^xQ
	(V_{u,n}^*\otimes G )
	=-V_{z,n}^T A^x
	\begin{pmatrix}
	V_{u,n}^*(1,:)\otimes G(1,:)\\
	\vdots\\
	V_{u,n}^*(3N,:)\otimes G(3N,:)
	\end{pmatrix},
	\end{aligned}
\end{equation}
where $ G=B^xV_{z,n} \in \mathbb{R}^{3N\times 3n} $ and the complexity of this operation is $ \mathcal{O}(Nn^3) $. Thus, the reduced matrix $ {Q}_{u,r} $ can be constructed without explicitly defining the matrix $ Q $. The transpose of the  Kronecker products of any given two vectors $ \mathbf{a} $ and $ \mathbf{b} $ can be represented as follows
	\begin{equation}\label{vec}
	\begin{aligned}
	(\text{vec}{(\mathbf{b}\mathbf{a}^\top)})^\top &=(\mathbf{a}\otimes \mathbf{b})^\top \\
	&=\mathbf{a}^\top\otimes \mathbf{b}^\top,
	\end{aligned}
	\end{equation}
	where vec $ (\cdot) $ denotes vectorization of a matrix. Using \eqref{vec}, the matrix $ N:=Q(V_{u,n}^*\otimes G )\in \mathbb{R}^{3N\times (3n)^2}$ in \eqref{goyal} can be constructed as follows
	$$
	N(i,:)=\left(\text{vec}\left(G(i,:)^T V_{u,n}^*(i,:)\right)\right)^T, \quad i=1,2,\ldots,3N.
	$$

In \cite{Benner19}, a pseudo-skeletal matrix decomposition  \cite{Mahoney09}, CUR, is used to increase further computational efficiency of the algorithm above. Instead of the CUR matrix decomposition, here we use "MULTIPROD" \cite{leva08mmm} in order to reduce
the complexity of the reduced nonlinear terms. MULTIPROD uses virtual array expansion to perform multiple matrix products. When the matrix $ V_{u,n}^*\in \mathbb{R}^{3N\times 3n} $ is reshaped as $ V_{u,n}^* \in \mathbb{R}^{3N\times 1 \times 3n} $,  then MULTIPROD is applied to $G$ and $V_{u,n}^*$ in $2$ and $3$ dimensions. MULTIPROD assigns virtually  a singleton to the third dimension of $G$, and we get the $3$-dimensional array (tensor) $\mathcal{N}:=\text{MULTIPROD}(G,V_{u,n}^*)\in\mathbb{R}^{3N\times 3n \times 3n}$. Thus, we can represent \eqref{goyal} as $Q_{u,r} = -V_{z,n}^T A^x N^{(1)}$, where $N^{(1)}\in \mathbb{R}^{3N\times (3n)^2}$ is the matricization of $\mathcal{N}$. In Section~\ref{num}, we compare the computational efficiency of computing the reduced matricized tensor $Q_{u,r}$ by the algorithm in \cite{Benner15}  with the algorithms in \cite{Benner18,Benner19} improved by the use of MULTIPROD.

\section{Numerical results}
\label{num}

We consider a test example for the SWE \eqref{swe1} on the spatial domain $\Omega = [0, 1]^2 $, with $ g = 1$, $f=0$, and with the initial conditions \cite{Sato18}
\begin{align*}
h({\mathbf x}, 0) &= 1+\frac{1}{2} \exp\left[-25\left(x-\frac{1}{2}\right)^2-25\left(y-\frac{1}{2}\right)^2\right],\\
u({\mathbf x}, 0) &= -\frac{1}{2\pi}\sin(\pi x) \sin(2\pi y),\\
v({\mathbf x}, 0) &= \frac{1}{2\pi}\sin(2\pi x) \sin(\pi y).
\end{align*}
The initial wave satisfy the periodic boundary conditions.
The final time is set to $T=50$, and spatial and temporal mesh sizes are taken as $\Delta x = 0.01$ and $\Delta t= 4 \Delta x$, respectively. This leads to a spatial grid of size $N=10000$, and $N_t=1250$ time intervals, so that each snapshot matrix $S_u$, $S_v$ and $S_h$ has size $10000\times 1250$.
For the FOM simulations, we consider the solutions of the SWE \eqref{semi-dis-poiss} in the Hamiltonian form with AVF time integrator (SWE-AVF), and the SWE \eqref{swef} in the f-plane with Kahan's time integrator (SWE-Kahan). In case of ROMs, related to the FOM \eqref{semi-dis-poiss} and by AVF time integrator, we consider the ROM \eqref{rom-semi} without DEIM approximation (POD-AVF) and the ROM \eqref{rom-semi-deim} with DEIM approximation (POD-AVF-DEIM). On the other hand, related to the FOM \eqref{swef} and by Kahan's time integrator, we consider the ROM \eqref{quad1} without tensorial framework (POD-Kahan) and with tensorial framework (TPOD-Kahan).
All the  simulations are performed on a machine with Intel$^{\circledR}$
Core$^{{\mathrm TM}}$ i7 2.5 GHz 64 bit CPU, 8 GB RAM, Windows 10, using 64 bit MatLab R2014.

\begin{figure}[htb!]
\centerline{\includegraphics[width=342pt,height=9pc]{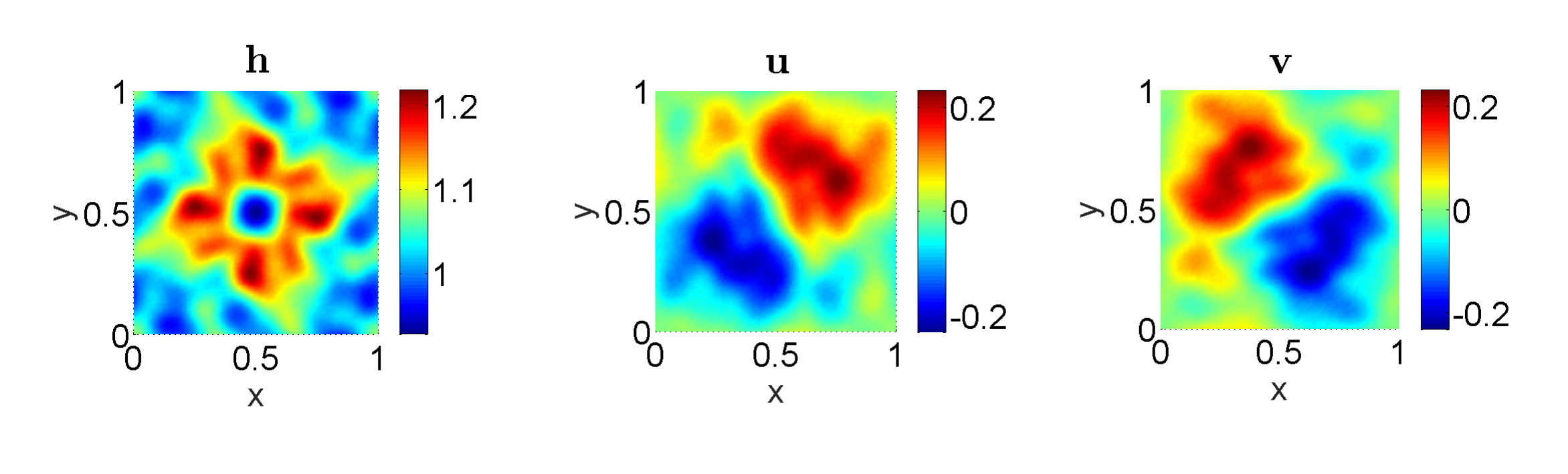}}
\centerline{\includegraphics[width=342pt,height=9pc]{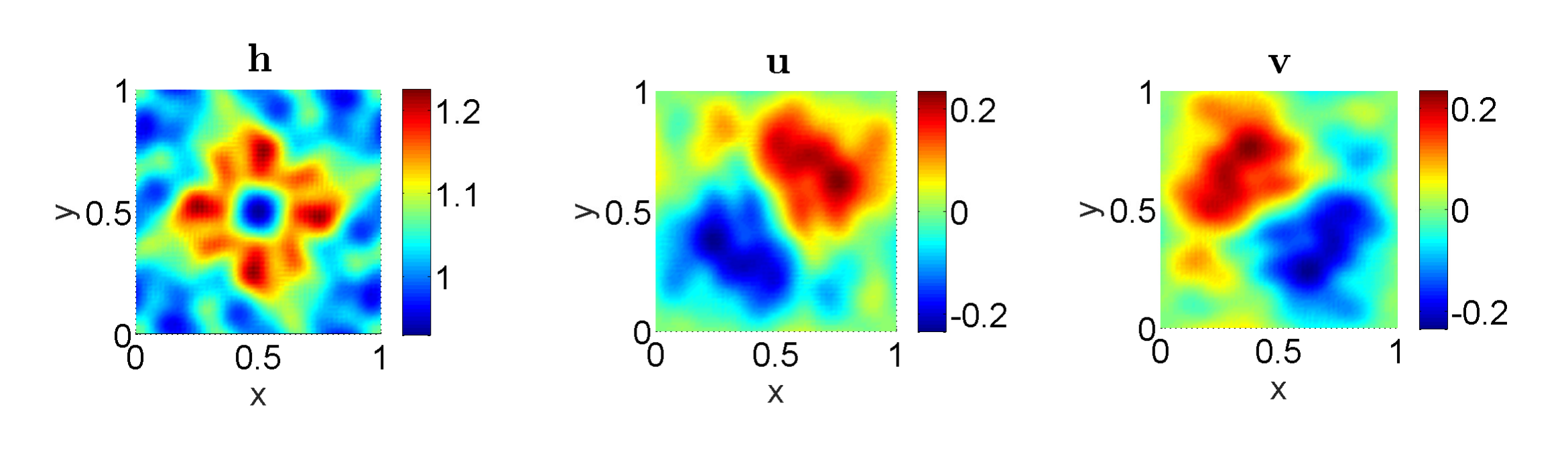}}
\caption{FOM solutions at final time $T=50$: (top) SWE-AVF, (bottom) SWE-Kahan.\label{foms}}
\end{figure}

\begin{figure}[htb!]
\centerline{\includegraphics[width=190pt,height=9pc]{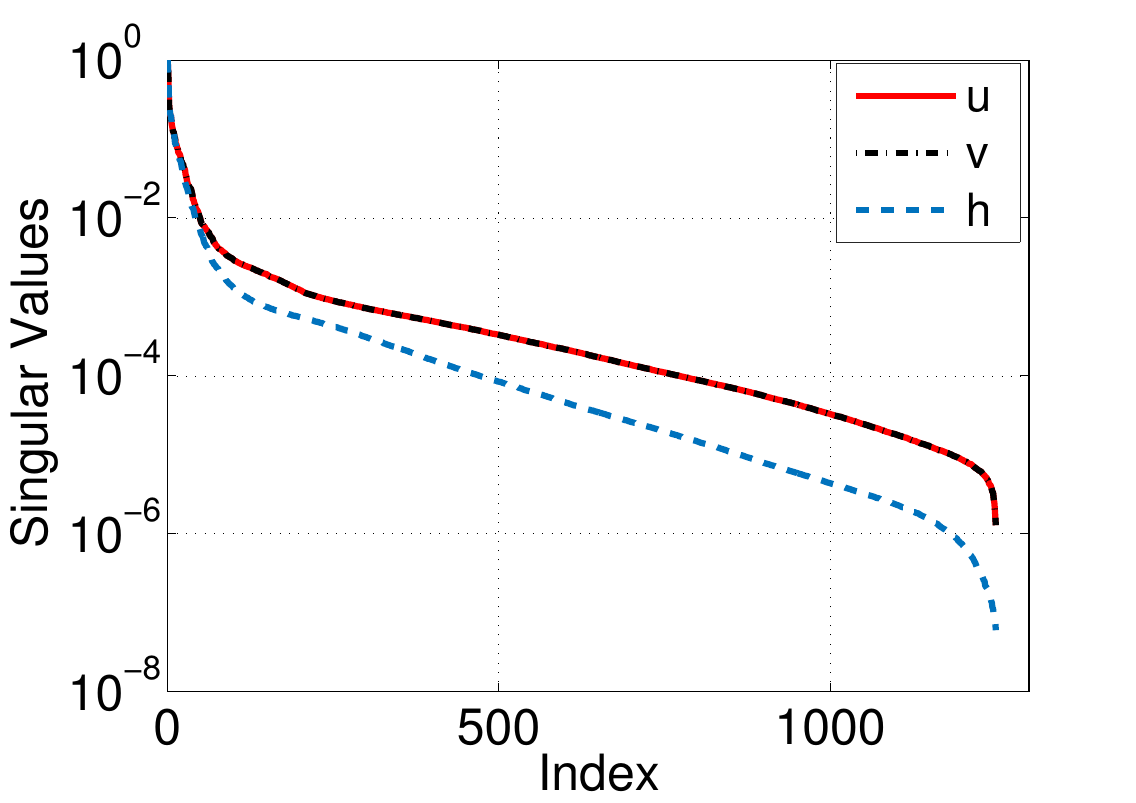}
\includegraphics[width=190pt,height=9pc]{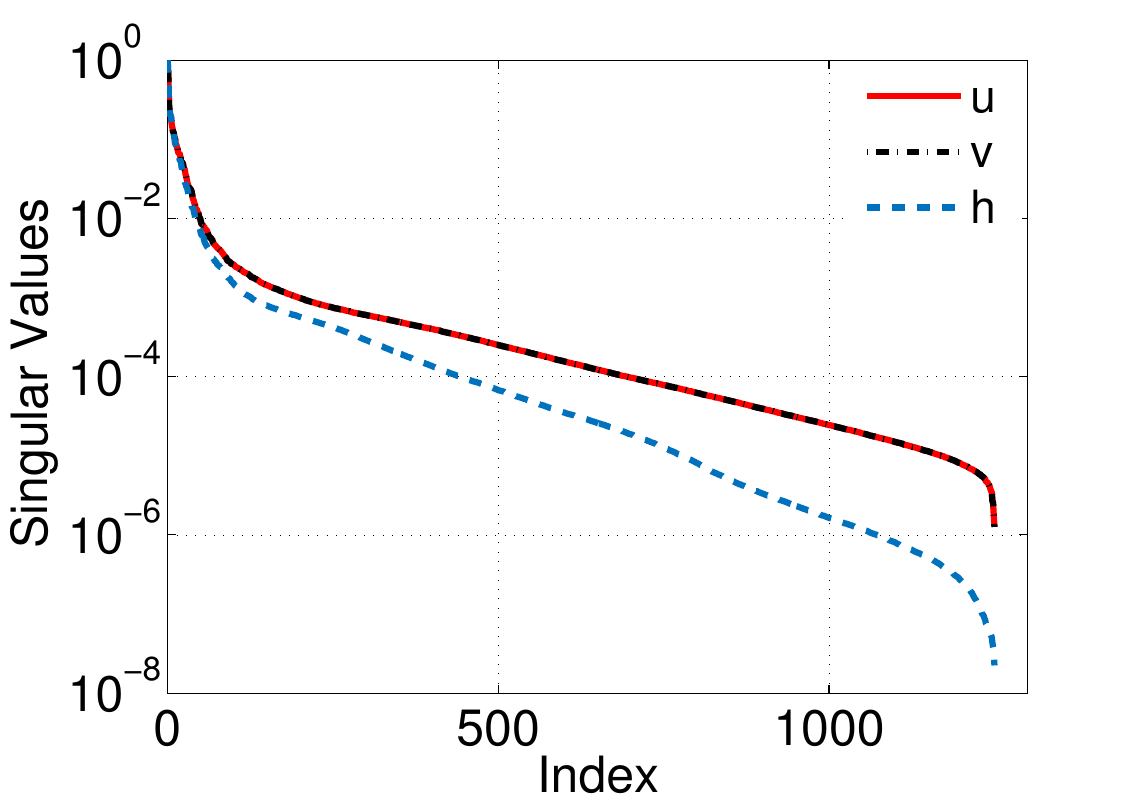}}
\caption{Singular values: (left) through SWE-AVF, (right) through SWE-Kahan.\label{svd}}
\end{figure}

The full order solutions of SWE-AVF and SWE-Kahan are given in Figure~\ref{foms}, depicting that the solutions by both approaches are in good agreement.
In Figure~\ref{svd}, we give the normalized singular values of the snapshot matrices $S_u$, $S_v$ and $S_h$ related to the velocity components $\mathbf{u}$, $\mathbf{v}$ and the height $\mathbf{h}$, respectively. The singular values decay slowly for each snapshot matrix obtained by both SWE-AVF and SWE-Kahan, which is the characteristic of the problems with wave phenomena in fluid dynamics \cite{Ohlberger16}. They pose a challenge for
reduced-order methods, since their dynamical behavior cannot be captured accurately
by the linear combination of a few POD modes.

In order to show conservation of the discrete energy \eqref{dener}, discrete enstrophy \eqref{denst}, and discrete vorticity \eqref{dmv} for a  FOM solution vector $\mathbf{z}$ (or ROM solution vector $\widehat{\mathbf{z}}$), we use the time-averaged absolute errors $\|\cdot\|_H$, $\|\cdot\|_Z$ and $\|\cdot\|_V$ defined by
\begin{align*}
\|\mathbf{z}\|_E = \frac{1}{N_t}\sum_{k=1}^{N_t}|E(\mathbf{z}^k)-E(\mathbf{z}^0)|, \quad E\equiv H,Z,V.
\end{align*}
On the other hand, to measure the global error between a discrete FOM solution vector and a discrete ROM approximation (FOM-ROM error), we define the following time averaged relative $L_2$ errors for the state variables ${\mathbf w} \equiv {\mathbf u},{\mathbf v},{\mathbf h}$
\begin{align*}
\|\mathbf{w}-\widehat{\mathbf w}\|_{Rel}=\frac{1}{N_t}\sum_{k=1}^{N_t}\frac{\|{\mathbf w}^k-\widehat{\mathbf w}^k\|_{L^2}}{\|{\mathbf w}^k\|_{L^2}}, \quad  \|{\mathbf w}^k\|_{L^2}=\sum_{i=1}^N{\mathbf w}^k_i\Delta x\Delta y,
\end{align*}
and to plot the FOM-ROM error at a specific time instance, we simply take the node-wise difference between the full discrete FOM solution and ROM approximation.

\begin{figure}[htb!]
\centerline{\includegraphics[width=190pt,height=9pc]{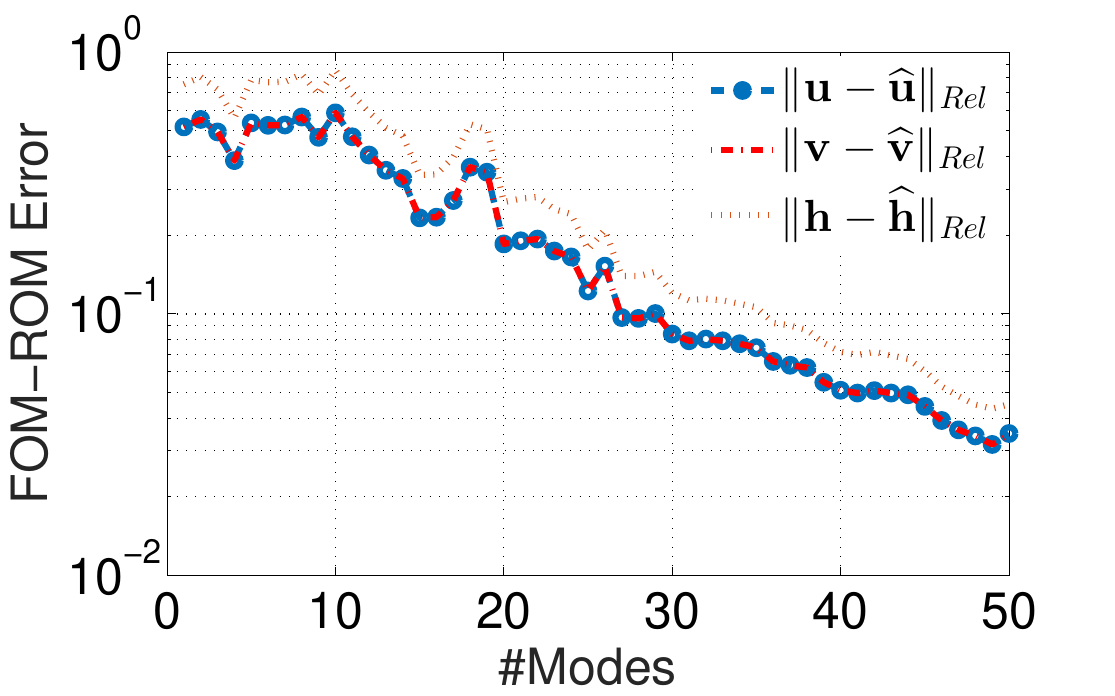}
\includegraphics[width=190pt,height=9pc]{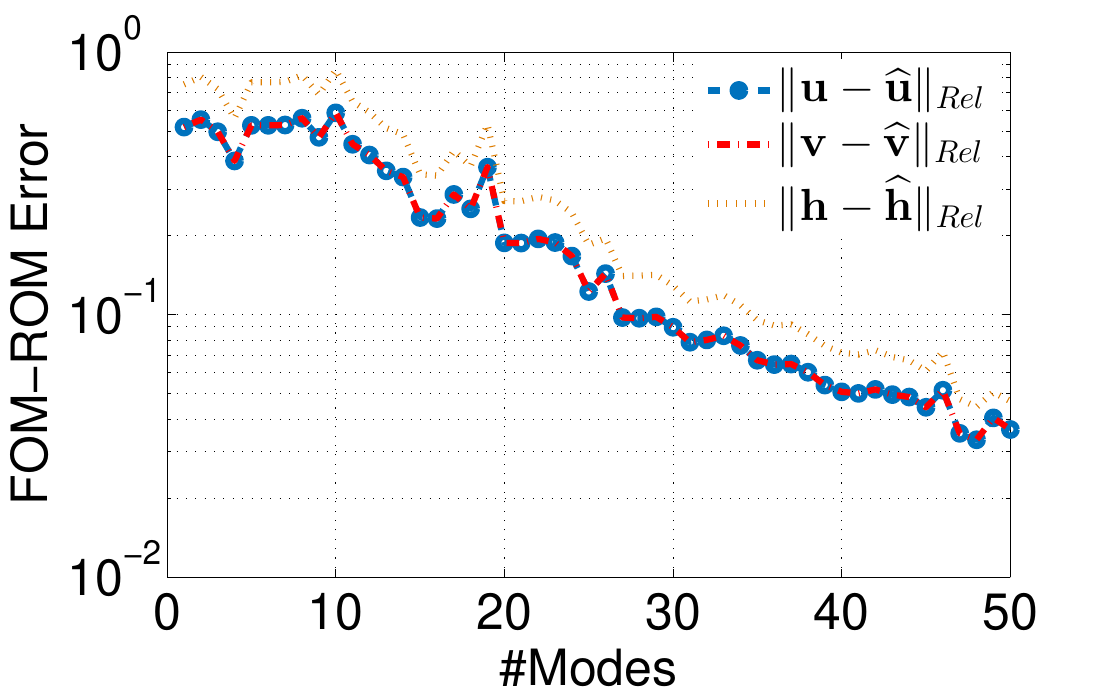}}
\caption{Time averaged relative $L_2$-errors vs. number of POD modes: (left) integration with AVF, (right) integration with Kahan's method. \label{l2l2}}
\end{figure}

Due to the slow decay of the singular values,  FOM-ROM errors for all state variable with varying number of POD modes in Figure~\ref{l2l2} decrease with small oscillations.  Therefore, for all the simulations, the number of POD modes is set to a relatively large value, $n=50$, according to the relative energy criteria \eqref{energy_criteria} with $\kappa =10^{-4}$. Using the same energy criteria, the number of DEIM modes is set to $m=90$.

\begin{figure}[htb!]
\centerline{\includegraphics[width=342pt,height=9pc]{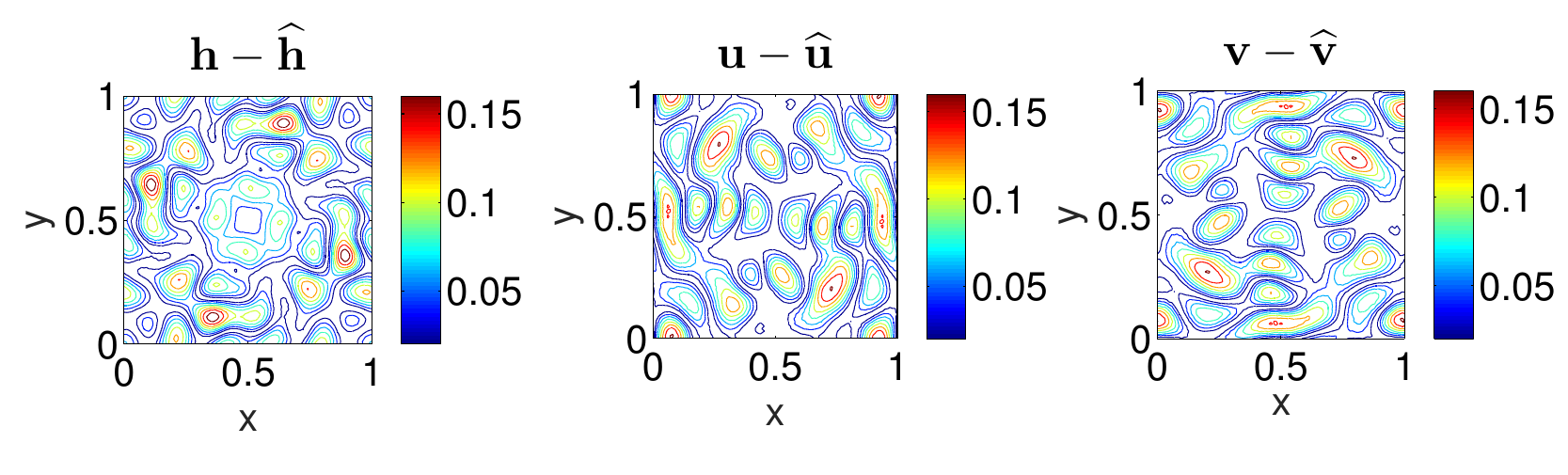}}
\centerline{\includegraphics[width=342pt,height=9pc]{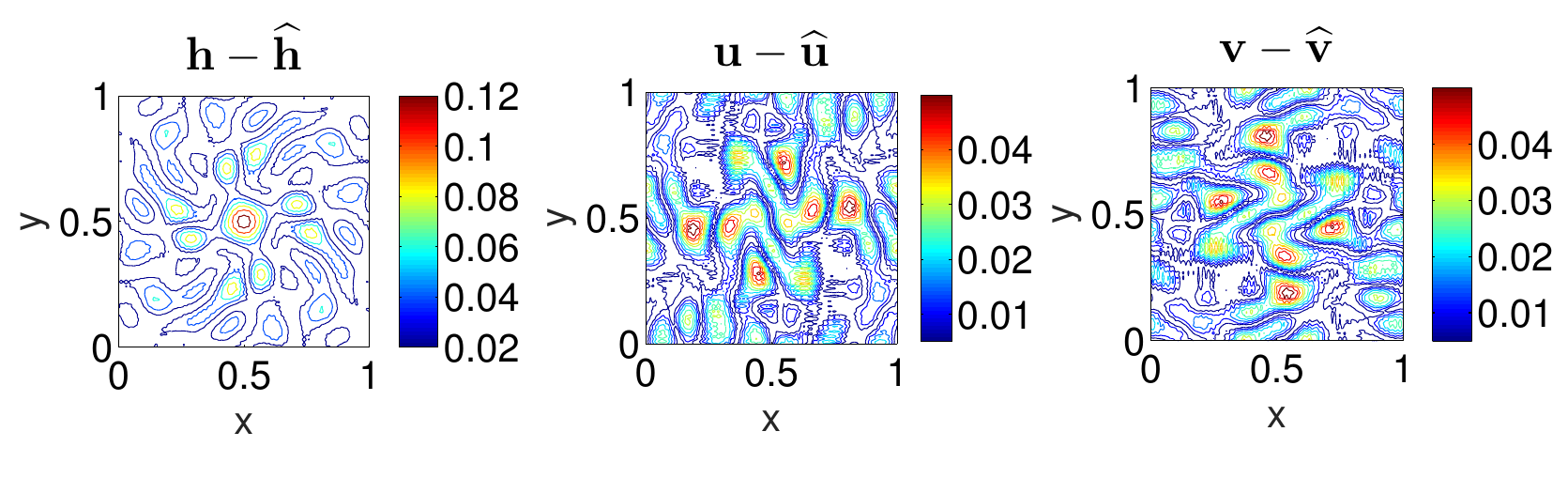}}
\caption{FOM-ROM errors at the final time $T=50$: (top) POD-DEIM-AVF, (bottom) TPOD-Kahan.\label{romerrors}}
\end{figure}

In Figure~\ref{romerrors}, the FOM-ROM errors at the final time for the height are at the same level of accuracy for both approaches, whereas the FOM-ROM errors of the velocity components by the TPOD-Kahan are smaller than the ones by the POD-DEIM-AVF, which might be due to differences in the solution algorithms.

\begin{figure}[htb!]
\centerline{\includegraphics[width=180pt,height=9pc]{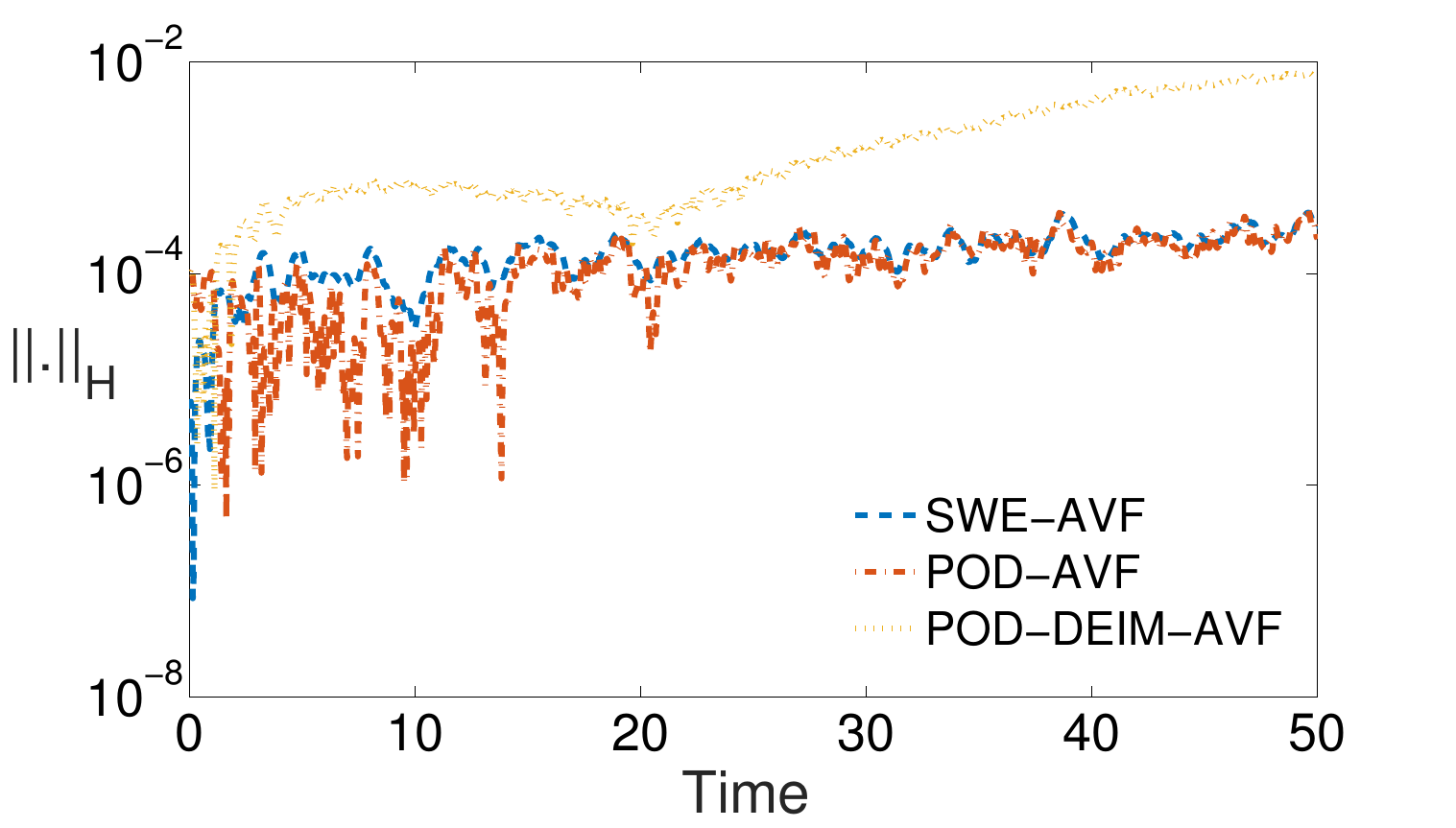}
\includegraphics[width=180pt,height=9pc]{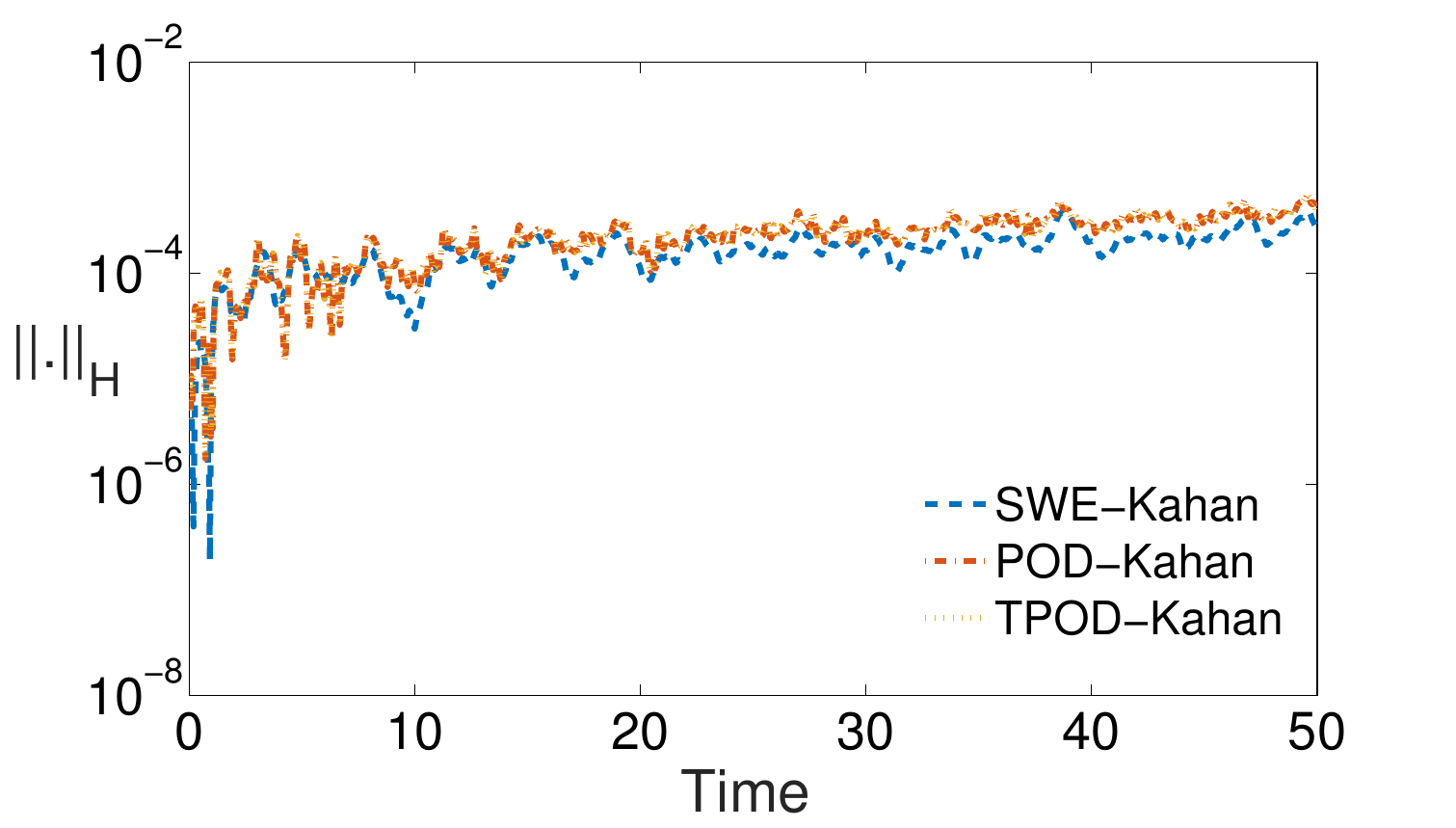}}
\caption{Energy error:  (left) integration with AVF, (right) integration with   Kahan's method. \label{enerror}}
\end{figure}

\begin{figure}[htb!]
\centerline{\includegraphics[width=180pt,height=9pc]{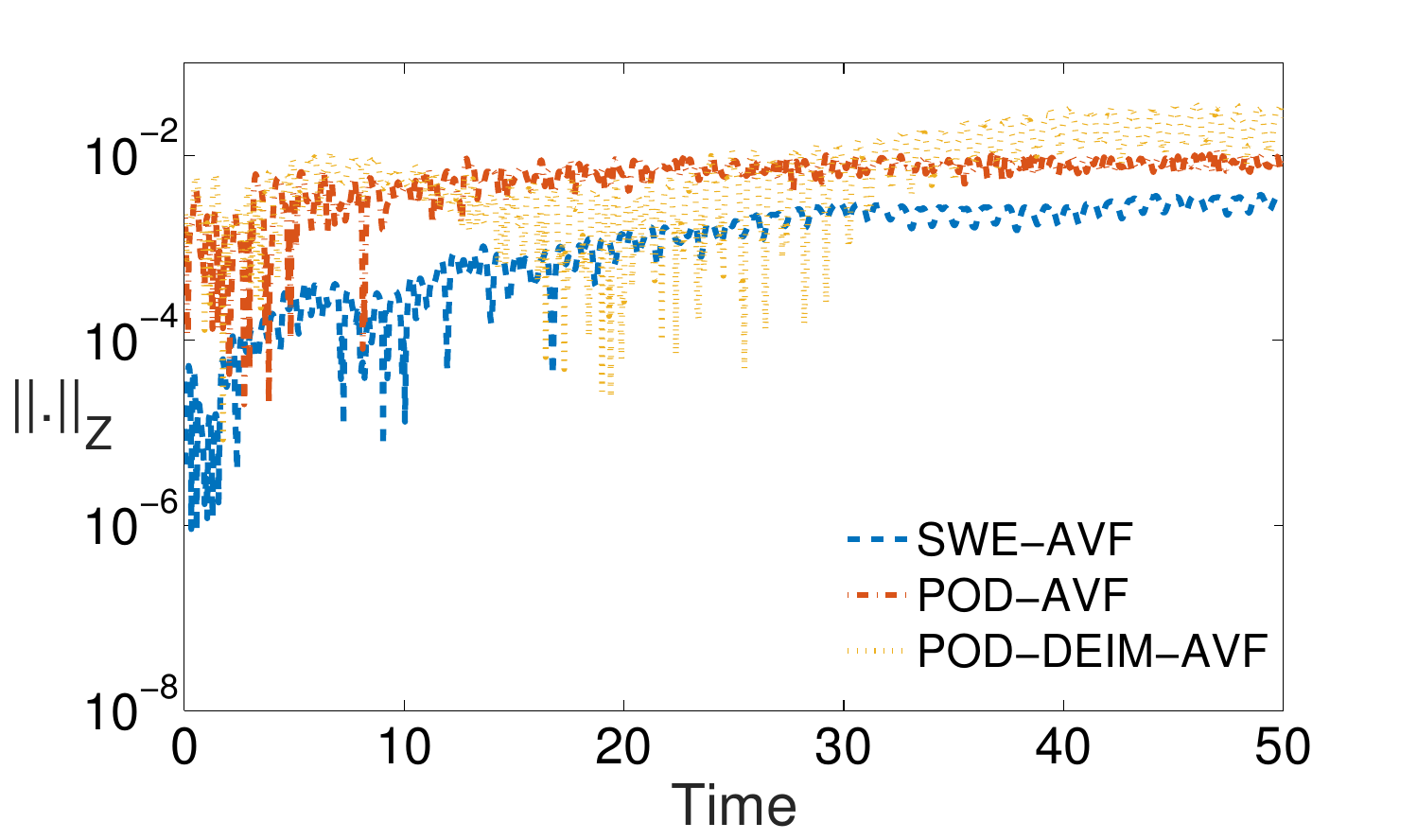}
\includegraphics[width=180pt,height=9pc]{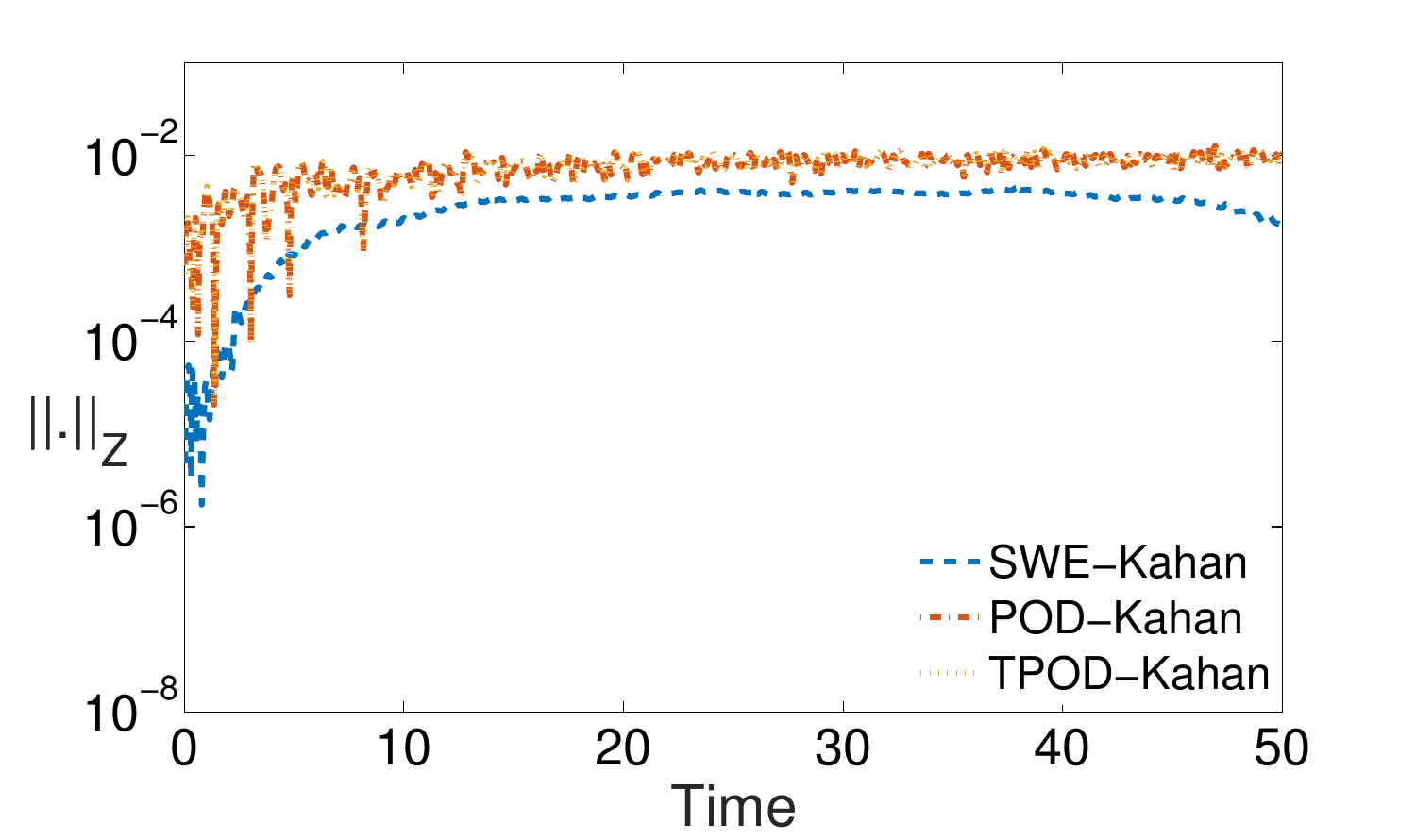}}
\caption{Enstrophy error: (left) integration
with AVF, (right) integration with Kahan's  method.  \label{enserror}}
\end{figure}

\begin{figure}[htb!]
\centerline{\includegraphics[width=180pt,height=9pc]{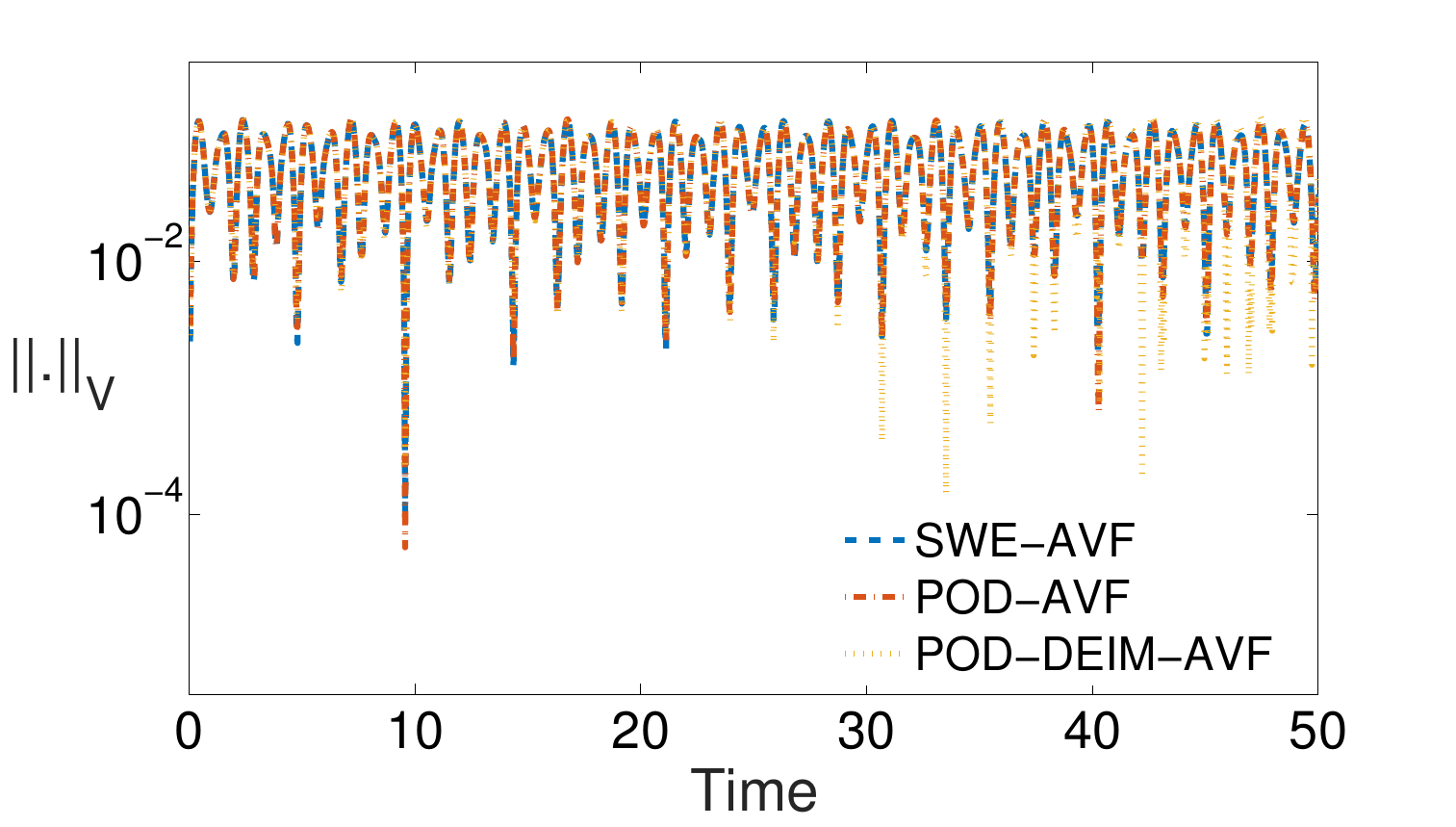}
\includegraphics[width=180pt,height=9pc]{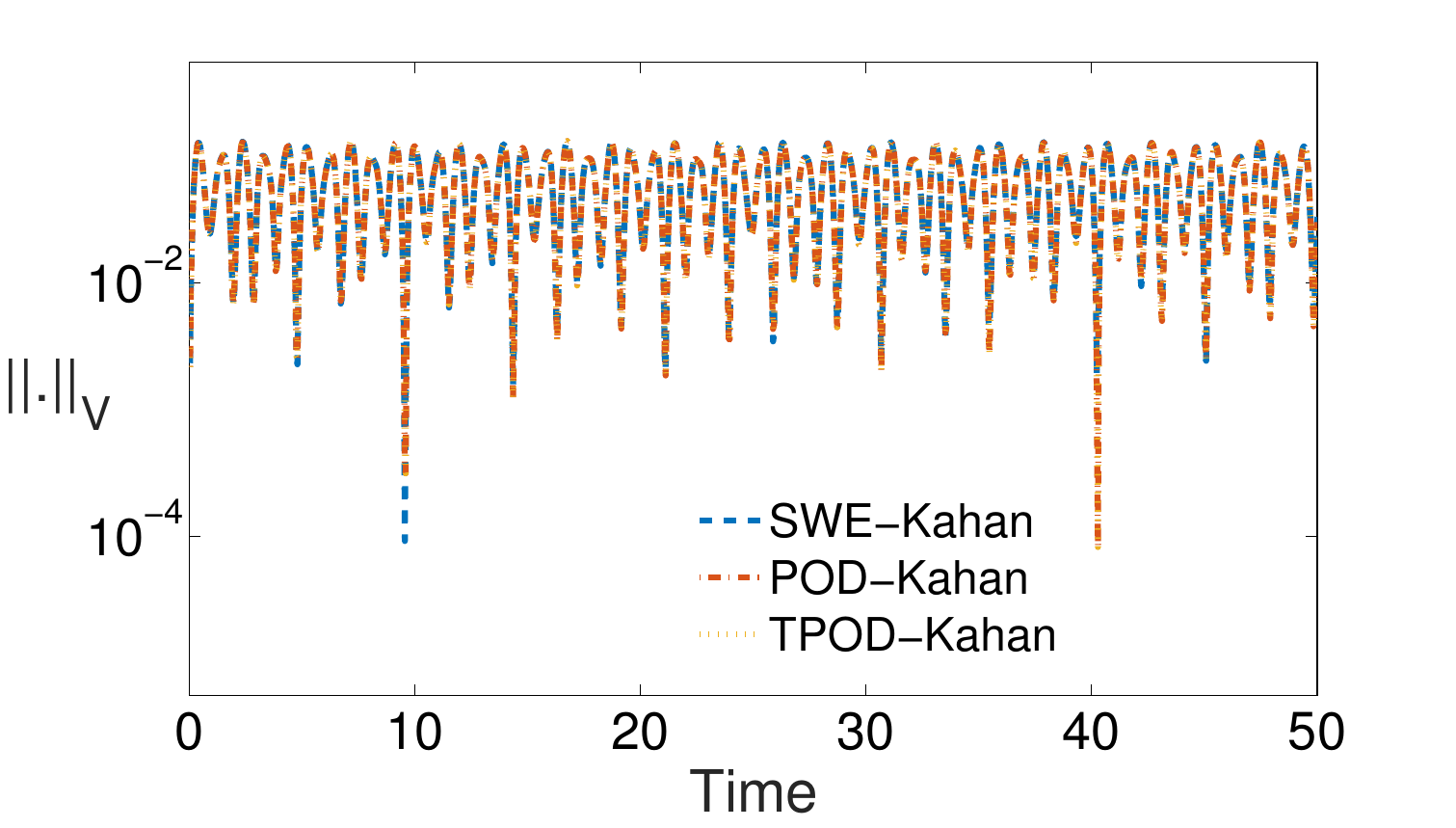}}
\caption{Vorticity error:  (left) integration
with AVF, (right) integration with Kahan's  method.     \label{vorerror}}
\end{figure}

In Figures~\ref{enerror}-\ref{vorerror}, the discrete energy error $|H^k -H^0|$, the discrete enstrophy error $|Z^k - Z^0|$, and the discrete vorticity error $|V^k -V^0|$ are plotted, $k=1,\ldots,N_t$.
The discrete energy, enstrophy, and vorticity are well preserved over the time interval. Because the vorticity is a quadratic conserved quantity, it is well preserved by both methods in Figure~\ref{vorerror}. The energy and the enstrophy errors in Figures~\ref{enerror}-\ref{enserror} show small drifts for POD-DEIM-AVF, but all they have bounded oscillations over the time, i.e., they are preserved approximately at the same
level of accuracy. The discrete energy and the discrete enstrophy is better preserved by TPOD-Kahan, since the reduced nonlinearity is not approximated  by  hyper-reduction. The mass is preserved up to machine precision since it is a linear conserved quantity, and it is not shown here.

\begin{table}[htb!]
\centering
\caption{Time averaged relative $L^2$-errors\label{solerr}}%
\begin{tabular}{llccc}
\hline
 &   & $\|\mathbf{u}-\widehat{\mathbf u}\|_{Rel}$  & $\|\mathbf{v}-\widehat{\mathbf v}\|_{Rel}$  & $\|\mathbf{h}-\widehat{\mathbf h}\|_{Rel}$ \\
\hline
POD-AVF & 30 POD modes  &  1.192e-01 & 1.192e-01  & 1.473e-02   \\
POD-DEIM-AVF     &  30 POD \&  90 DEIM modes  & 2.874e-01  & 2.874e-01  & 3.478e-02   \\
\hline
TPOD-Kahan & 30 POD modes & 1.265e-01  &1.265e-01   &1.567e-02 \\
\hline
\end{tabular}
\end{table}

\begin{table}[htb!]
\centering
\caption{Time averaged absolute errors for the conserved quantities\label{conserr}}%
\begin{tabular}{llccc}
\hline
 &   & \textbf{Energy}  & \textbf{Enstrophy}  & \textbf{Vorticity} \\
\hline
POD-AVF  &  30 POD modes  &   3.768e-05 & 7.871e-03           & 2.123e-05   \\
 POD-DEIM-AVF       & 30 POD\&   90 DEIM modes & 1.968e-03 & 5.137e-03           & 1.042e-03   \\
\hline
TPOD-Kahan & 30 POD modes   &  2.901e-05         &   3.108e-03           &    3.454e-05   \\
\hline
\end{tabular}
\end{table}

\begin{table}[htb!]
\caption{CPU time (in seconds) and speed-up factors \label{cpu}}
\centering
\resizebox{15cm}{!}{%
\begin{tabular}{llllllll}
\hline
\multicolumn{4}{c}{\textbf{POD-DEIM-AVF}} & \multicolumn{4}{c}{\textbf{TPOD-Kahan}} \\
\hline
 &  & CPU time & speed-up & &  & CPU time & speed-up \\
\hline
FOM   &   & 529.80  &   & FOM  &    & 788.89 &   \\
\hline
 \multirow{3}{*}{POD}  & basis computation  & 33.46  &  & \multirow{3}{*}{POD}  &  basis computation  & 31.02 &   \\
	 & online computation  & 291.69  &  &  &  online computation & 289.01 &   \\
	& total  & 325.15  & 1.63 &   &  total  & 321.03 & 2.45  \\
\hline
 \multirow{3}{*}{DEIM}  & basis computation  & 33.89  & &  \multirow{3}{*}{TPOD}  &  tensor computation \cite{Benner15}(MP)  & 20.78 (9.43) &   \\
	 & online computation  & 24.10  &  &  &  online computation & 6.42 &  \\
	& total  & 57.99  & 9.13 & &  total \cite{Benner15}(MP)  & 27.20 (15.89) &  29.00 (49.64) \\
\hline
\end{tabular}}
\end{table}

The time-averaged relative $L_2$-errors between FOM and ROM solutions in Table~\ref{solerr} are at the same level of accuracy for both versions of ROMs.  The ROMs for the height are more accurately resolved than for the velocity components.
The conserved quantities are more accurately preserved by the POD-AVF  and by the TPOD-Kahan than by the POD-DEIM-AVF in Table~\ref{conserr}. This indicates that using DEIM, the stability of the ROM solutions can not be guaranteed in long term integration, which is reflected in the preservation of the conserved quantities.

The CPU times and speedup factors in Table~\ref{cpu} for $n=50$ POD and $m=90$ DEIM modes, show the computational efficiency of the TPOD-Kahan over POD-DEIM-AVF due to the separation of offline-online computation. In Table~\ref{cpu}, basis computation includes SVD computation, and online computation consists of time needed for projection and solution of reduced system. As pointed out in \cite{Kramer19} that for some problems, a large number of DEIM interpolation points are required to achieve accurate solutions, which increases the computational cost of the ROMs in online computation, as in Table~\ref{cpu}. The computational efficiency is further increased by exploiting the sparse matrix structure of the discretized SWE using MULTIPROD (MP) in the algorithm of \cite{Benner18,Benner19} over \cite{Benner15} as shown in Table ~\ref{cpu}.

\begin{figure}[htb!]
\centerline{\includegraphics[width=380pt,height=11pc]{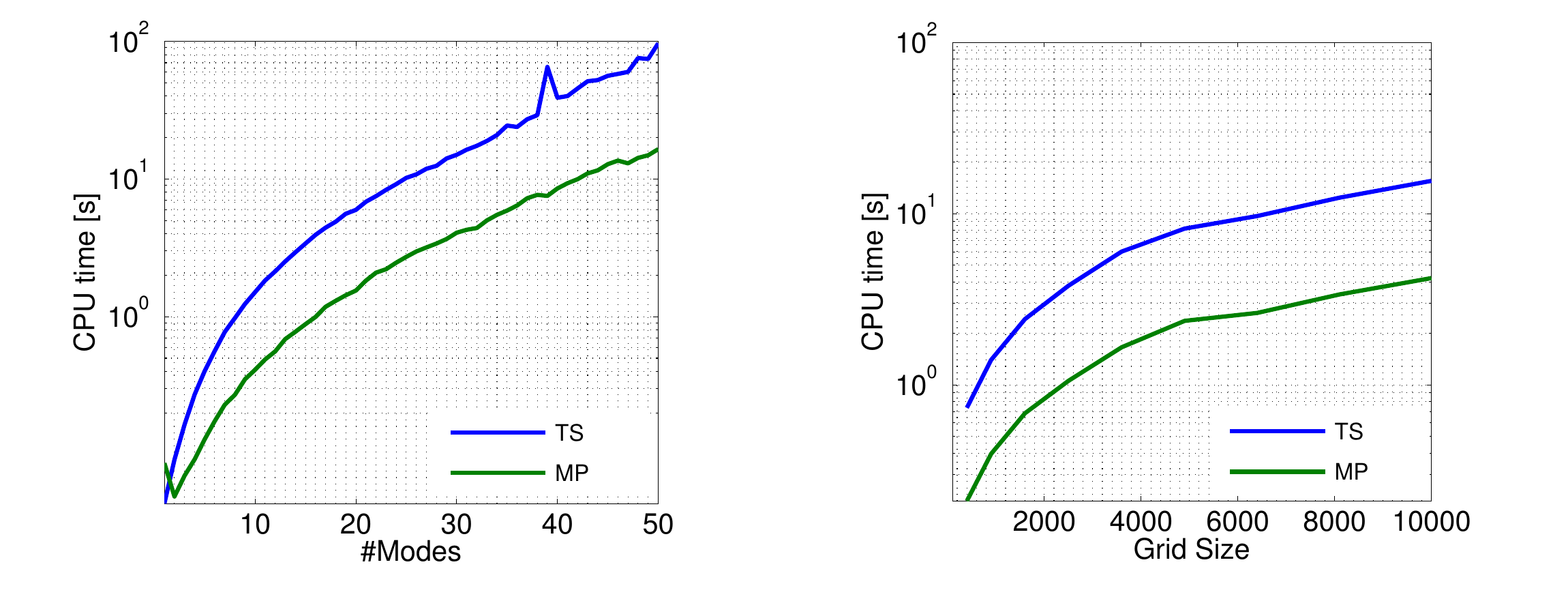}}
\caption{Tensor calculations by  \cite{Benner15} (TS) and MULTIPROD (MP): (left) FOM is fixed with $N=10000$, (right) ROM is fixed with $n=50$.\label{bg}}
\end{figure}

For the reduced matricized tensor $Q_{u,r}=-V_{z,n}^TA^xQ\left(V_{u,n}^*\otimes (B^xV_{z,n}) \right)$, the computational time by the two-sided projection method in \cite{Benner15} (TS) is compared with the algorithms in \cite{Benner18,Benner19} using MULTIPROD. In Figure~\ref{bg}, left, we give the required computational time versus the
number of POD modes by fixing the number of grid points $N=10000$, whereas the required computational time for varying number of grid points using a fixed number of POD mode $n=50$ are given in Figure~\ref{bg}, right. Both figures show the computational efficiency using MULTIPROD by increasing the size of the FOM and ROMs.

\section{Conclusions}
\label{conc}

We have shown that the two kind of ROMs preserve the conserved quantities of the SWE and yield stable reduced solutions in long time integration. For both approaches, the reduced  solutions and the reduced conserved quantities have almost the same level of accuracy as the full ones.  Using TPOD-Kahan and exploiting the quadratic structure of the SWE, the online computation time of ROMs can be much reduced compared to the POD-DEIM-AVF. The online computation time of the ROMs is further reduced by applying the algorithms in \cite{Benner18,Benner19} in combination with the MULTIPROD \cite{leva08mmm}. In a future work, we will investigate both approaches for the SWE with full Coriolis force \cite{Stewart16} and thermal SWE \cite{Eldred19}.

\subsubsection*{ACKNOWLEDGEMENT} The authors thank for the constructive comments of the referees which helped much to improve the paper.  This work  has been supported by 100/2000 Ph.D. Scholarship Program of the Turkish Higher Education Council.


\begin{tcolorbox}
{\bf How to cite this article} Karasözen B, Yıldız S, Uzunca M. Structure preserving model order reduction of
shallow water equations. {\it Math Meth Appl Sci}. 2020; 1-17. \href{https://doi.org/10.1002/mma.6751}{https://doi.org/10.1002/mma.6751}
\end{tcolorbox}
 
\end{document}